
\input amstex 
\documentstyle{amsppt}
\magnification 1100
\TagsOnRight
\NoBlackBoxes
\hsize=6.30 truein
\hcorrection{.120in}
\vsize=8.5 truein

\def \RR   {{R\!\!\!\!\!I~}} 

\def\Box{\sqcup\!\!\!\!\sqcap}

\noindent{\bf Published as :}  
       P.G. LeFloch, Propagating phase boundaries: formulation of the problem and existence via the Glimm scheme, 
      {\it Arch. Rational Mech. Anal.} {\bf 123} (1993), 153--197

      \
      
      \
      
{\baselineskip.4cm
\centerline{\bf PROPAGATING PHASE BOUNDARIES:}
\vskip2mm
\centerline{\bf FORMULATION OF THE PROBLEM AND EXISTENCE}
\vskip2mm
\centerline{\bf VIA GLIMM METHOD}
\vskip.4cm 
\vskip.4cm
\centerline{Philippe G. LeFLOCH}
\vskip.2cm
\centerline{\bf Address from 1990 to 1992 :} 
\centerline{Courant Institute of Mathematical Sciences}
\centerline{New York University}
\centerline{251 Mercer Street, New York, NY 10012}
\vskip1mm

\


\

\centerline{\bf Current address : } 
\centerline{Laboratoire Jacques-Louis Lions}

\centerline{Centre National de la Recherche Scientifique}

\centerline{Universit\'e de Paris 6}

\centerline{4 place Jussieu, 75252 Paris, France.}

\centerline{E-mail: {\it LeFloch\@ann.jussieu.fr}} 

\vskip.5cm
\noindent{\bf Abstract. }
In this paper we consider the hyperbolic-elliptic system of two conservation
laws that describes the dynamics of an elastic material governed by a
non-monotone strain-stress function. Following Abeyaratne and Knowles, 
we propose a notion of admissible weak solution for this system
in the class of functions with bounded variation. The formulation includes
an entropy inequality, a kinetic relation (imposed along any subsonic
phase boundary), and an nucleation criterion (for the appearance of
new phase boundaries).  We prove the $L^1$--continuous dependence of the
solution to the Riemann problem. Our main result yields the existence
and the stability of propagating phase boundaries. The proof is
based on Glimm's scheme and in particular on techniques going back to 
Glimm-Lax. In order to deal with the kinetic relation, we prove 
a result of pointwise convergence for the phase boundary.
\vfill
%

















\

\noindent{\bf 0. INTRODUCTION}
\vskip.2cm

This paper deals with the following system of two conservation laws which
describes the motion of an elastic material
$$
\partial _t w -\partial _x v = 0, \qquad 
\partial_t v -\partial _x \sigma(w)=0.
\leqno{(0.1)}
$$
Here $w > -1$ and $v$ represent the deformation gradient and the velocity
of the material, respectively. The stress function $\sigma : \, ]-1,\infty[
\to \RR$ is assumed to be monotonically increasing except in an interval
$]w_M,w_m[$.
Such a form of the stress function is
typical in the modeling of solid materials which admit different phases.
A van~der~Waals gas also is described by a very similar system. System (0.1)
is of mixed type, i.e. {\it hyperbolic} in the phase~1 region
${\Cal  H}_1 = \{ w < w_M\}$ and in the phase~3 region ${\Cal H}_3 =\{
w>w_m\}$, but {\it elliptic} in the intermediate region of phase 2~states.
The phase 2~states are known to be both mathematically and physically
unstable (James [23]). We will consider here exclusively solutions which
take their values in the {\it phase~1} or {\it phase~3} regions {\it only}.

Solutions to (0.1) in general are discontinuous and so must be understood
in the sense of distributions; see Lax [25], [26] for background on weak 
solutions. Such discontinuous solutions are in general non-unique
and those having a physical meaning must be selected through an admissibility
(or entropy) criterion. We refer to Dafermos [7] for a review of entropy conditions
in the setting of hyperbolic problems. As was pointed out by James
 [23], the mixed system (0.1) possesses
a high degree of {\it non-uniqueness}, that a number of authors
have attempted to resolve by means of suitable generalizations
of entropy criteria from the theory of hyperbolic conservation laws.
First of all, Shearer [38] considered the Lax entropy criterion [25], [26].  The viscosity 
and viscosity-capillarity approches have been analyzed by Slemrod 
[41], [42]; cf.~also Hagan-Slemrod [18], Pego [35] and Shearer [39], [40].  
Next, Hattori [19], [20] has investigated the application to (0.1)
of the entropy rate admissibility criterion proposed by Dafermos [6].
Hsiao [22] has considered the Liu entropy criterion [32] which allows one to treat
equations of state losing genuine nonlinearity in 
hyperbolic regions. Another approach to resolve
the non-uniqueness can be found in a work by Keyfitz [24]. Additional 
material on system (0.1) is found in [12], [13] and [36].

All the above works consider the Riemann problem only, i.e. a Cauchy
problem for (0.1) with initial condition which consists of two constant states.
This problem can be solved explicitly (in a possibly non-unique way) by using simple
waves (i.e. shock waves, rarefaction waves or contact discontinuities).
Adding an ``admissibility criterion'' allows one to reduce the class
of (admissible) solutions and in most situations to select a unique
solution.  However, it must be emphasized that the solution
of the Riemann problem (when it is unique) depends on the
chosen admissibility criterion. It turns out that there is no preferred
criterion for the selection of the ``physically meaningful'' solutions
of (0.1).

A different approach was recently investigated by Abeyratne and Knowles
in [2]. The main suggestion of these authors is that system (0.1)
is not physically complete enough to describe the evolution of a phase
boundary in an elastic material. It must be completed with a {\it kinetic
relation} imposed along any subsonic phase boundary: this kinetic
relation actually yields the rate of entropy dissipation across the
phase discontinuity. Moreover, Abeyaratne and Knowles add an {\it initiation
criterion} which controls the possible appearance of a new phase. We refer
to [1] and the references therein for the motivation of introducing
a kinetic relation and an initiation criterion which are actually classical in the context of
quasi-static problems. Cf. also Gurtin [17] and Truskinovsky [44] for related ideas.

Abeyaratne and Knowles
proved in [2] that the Riemann problem for (0.1) always admits a unique
admissible solution, i.e. a weak solution satisfying the kinetic relation
and the initiation criterion, as well as the entropy inequality which 
reads
$$
\partial_t \big(W(w) +{v^2\over 2}\big) -\partial_x \big(\sigma(w)v\big)\le 0,
\leqno{(0.2)}
$$
where $W:\, ]-1,\infty[\to\RR$ is the internal energy function defined by
$$
W(w) = \int_0^w \sigma(y)\, dy \quad\hbox{ for }\quad w\in ]-1,\infty[.
\leqno{(0.3)}
$$
Next they showed in [3] that the solution of the Riemann problem found
by Slemrod through the viscosity-capillarity approximation corresponds to a 
special choice of kinetic relation in their approach.
It is not difficult to check also that the solution found by Shearer
 using Lax entropy inequalities coincides with the maximally dissipative
kinetic relation investigated in [4]. (I thank Michael Shearer 
for pointing that out to me.)

The present paper is devoted to continuing the analysis of system (0.1)
through the approach of Abeyaratne and Knowles. As in [2], we will
restrict ourselves to the case of a piecewise linear stress-function.
This assumption simplifies the calculations but it is not a real
restriction to the results of this paper.

Our purpose is first (Sections~1 and 2) to give a slightly different
presentation of the ideas of [2], which as we think clarifies the concepts
of kinetic relation and initiation criterion introduced by Abeyaratne
and Knowles. Section~1 presents the mathematical formulation of a
well-posed (at least for Riemann data) problem associated with system
(0.1). As is usual for hyperbolic problems, we consider bounded solutions of
bounded variation (BV).  Our formulation follows [2] with
however two main modifications. The kinetic relation is introduced
from a completely {\it dynamical point of view} and not as a
generalization of the quasi-static point of view as was done in [2].
That leads us to a {\it larger range} of admissible values for the  --- as we 
call it below --- entropy dissipation function in the kinetic relation. 
Furthermore, the initiation criterion at some point $x$ is formulated here 
in two different ways, depending on whether $x$ is in an interior point of
the space interval $[a,b]$ where we set the problem, or $x$ is a point  
of its boundary. For definiteness, we allow {\it spontaneous  initiation of a new phase} 
only at the extremities of the bar $[a,b]$, which is consistent with
the classical static theory.

Then Section~2 describes briefly the solution of the Riemann problem.
We explain how to take into account the two changes above in the
construction of [2]. The main result of this section establishes
the $L^1$ {\it continuous dependence} of the Riemann solution with 
respect to its initial states.  It must be emphasized that the two 
observations above are essential for the continuous dependence property 
to hold, especially our condition that a new phase may occur spontaneously
only at the extremities of the bar. The same results are also obtained
for the Riemann problem in a half-space.
Note that, although uniqueness of the admissible solution holds
for the Riemann problem, nothing is known for the general Cauchy
problem. As a matter of fact, the issue of uniqueness for conservation 
laws is understood in a few number of situations only. (See, for hyperbolic
 problems, LeFloch-Xin [31] and the references therein.)

The second part of the paper (Sections~3 and 4) focuses on the solutions
of the Cauchy problem for system (0.1), which are BV perturbations
of a single propagating phase boundary separating a phase~1 state and a
phase~3 state. We prove the existence of admissible weak solutions of this
form, when the initial data on both sides of the phase discontinuity
has small total variation. We treat the case of any non-characteristic phase
boundary as well as the case of a characteristic phase boundary
provided that no strong wave arises from perturbating the states
on both sides of the phase boundary. The random-choice scheme due to Glimm [15] is used
to construct approximate solutions to the problem. Its stability in the 
BV norm is proved from an essentially  linear estimate of wave interactions between
two Riemann solutions. Such linear interaction terms were 
used in a different situation by Chern [5] and Schochet [37]. Note that 
the strength of the phase discontinuity is not (and can not be) assumed to be small in any sense.

The stability of the scheme in the total variation norm is sufficient
to extract a subsequence converging to a weak solution of the
problem. This convergence result holds almost everywhere with respect to 
the Lebesgue measure. This is sufficient to show that the scheme converges 
to a weak solution of the problem. But, proving that this solution is
 {\it admissible} requires a result
of {\it pointwise convergence} of the phase boundary. In Section~4, we 
establish this property by using the technique of analysis due
to Glimm-Lax [16]. We next prove that it is sufficient, at least for 
non-stationary phase boundaries, for the passage to the limit in the
 kinetic relation.

An extension of the results in this paper to arbitrary large initial 
data would require a better understanding of the phenomena of initiation of new phases. 

Many ideas in this paper are related to those in the developing theory of
nonlinear hyperbolic systems in non-conservative form for which we
refer the reader to Dal~Maso--LeFloch--Murat [8] and LeFloch--Liu [30]; see also [27] to [29].

\

\noindent{\bf Acknowledgements.} This work was done while the author was 
assistant professor at CIMS (Courant Instructor), and on leave from a CNRS research
position at the Ecole Polytechnique (France). This work was in part supported
 by a NSF grant DMS-88-06731 and by the CNRS (Centre National de 
la Recherche Scientifique). It is a great pleasure for me to thank 
Bob Kohn who encouraged me to
study this problem and accepted the job of reading the first version 
of this paper. The hospitality of Peter Lax at Courant Institute
was highly appreciated. I am also graceful to Rohan Abeyaratne and Michael
 Shearer for discussions.

\

\




%
\noindent{\bf 1. MATHEMATICAL FORMULATION OF THE PROBLEM}
\vskip.2cm

This section describes the formulation of the Cauchy problem associated with
the mixed system (0.1). The formulation includes the system of conservation
laws (mass, momentum) (0.1) together with the (Clausius-Duhem)
entropy inequality associated
with the entropy $W(w)+{v^2\over 2}$. It is made complete by adding
to these both a kinetic relation along any subsonic phase boundary
and an initiation criterion for the occurrence of possible new phase
boundaries in the solution. We specify below the assumptions on
the kinetic relation and the initiation criterion which will be essential
to the results of Section~2.  This section also introduces notation
which will be of constant use throughout this paper.

We write system (0.1) in the form
$$
\partial_t u +\partial_x f(u) = 0, \quad
u={v\choose w},  \quad f(u) = {-\sigma(w)\choose -v}.
\leqno{(1.1)}
$$
For simplicity, we shall assume that the stress-function $\sigma\colon ]-1,\infty[
\to\RR$ is a piecewise linear function of the following form 
$$
\sigma(w) = \cases
k_1w  \hfil & \text{ for } \quad -1\le w \le w_{M}, 
  \cr
 k_3 w_m + (k_1w_M-k_3 w_m) (w-w_m) /(w_M-w_m)
\hfil & \text{ for } \quad w_M \le w\le w_m, 
   \cr
k_3 w \hfil & \text{ for } \quad w_m\le w. \cr
\endcases
\leqno{(1.2)}
$$
The constants $k_1,k_3,w_m$ and $w_M$ in (1.2) are assumed to satisfy
the properties
$$
0 < k_3 < k_1 \quad\hbox{ and }\quad 0 < w_M < w_m.
\leqno{(1.3)}
$$
We shall use the notation
$$
\sigma_M = k_1 w_M \quad\hbox{ and }\quad \sigma_m = k_3 w_m.
$$
The phase~1 region ${\Cal H}_1 = \{-1 < w \le w_M\}$ and the phase~3 region
${\Cal H}_3 = \{w\ge w_m\}$ correspond to observable and stable states.
In our formulation below, the solution {\it cannot enter} the unstable phase~2
region $\{w_M < w < w_m\}$ and so must jump from ${\Cal H}_1$ to ${\Cal H}_3$
or conversely. A discontinuity between two states in different phases
is called a phase boundary.

\



System (1.1) is linear hyperbolic in ${\Cal H}_1$ and ${\Cal H}_3$, 
and the corresponding characteristic speeds are
$$
\pm c_1 =\pm \sqrt{k_1} \,\,\hbox{ in }\,\, {\Cal H}_1
\quad\hbox{ and }\quad \pm c_3 = \pm \sqrt{k_3} \,\,\hbox{ in }\,\, 
{\Cal H}_3.
$$
In view of (1.3), the waves in phase~1 travel faster than those in phase~3,
i.e. $c_1 > c_3$.
We may also use the notation
$$
c(w) =c_1 \quad\hbox{ if }\quad w \le w_M, \quad c_3\quad \hbox{ if }\quad
w\ge w_m.
\leqno{(1.4)}
$$
Note that $c(w)$ is not defined if $w$ belongs to $]w_M,w_m[$. With some
abuse of notation, ${\Cal H}_1$ and ${\Cal H}_3$ will sometimes also
denote $\{ (v,w) / -1 < w \le w_M, \, v\in \RR\}$ and $\{(v,w)/
w_m \le w, \, v\in \RR\}$, respectively. We also set ${\Cal H} = {\Cal H}_1
\cup {\Cal H}_3$.

Since (1.1) is linear hyperbolic in ${\Cal H}_1$ and ${\Cal H}_3$,
possible discontinuities in the initial data for (1.1) are simply
advected along the characteristic lines of slopes either $\pm c_1$
or $\pm c_3$. (This is true at least up to the time of appearance of
a new phase.) The elementary waves in each of the regions ${\Cal H}_1$ and
${\Cal H}_3$ are contact discontinuities.
Hence, the special choice (1.2) for the constitutive law is very convenient. 
It makes quite simple the analysis in the hyperbolic regions and allows
us to focus on the phase boundaries between ${\Cal H}_1$ and ${\Cal H}_3$.
We will see that the description of the appearance and the
evolution of the phase boundaries is far from trivial.

In the theory of hyperbolic conservation laws, it is standard to consider
solutions $u=(v,w)$ to (1.1) in the functional space 
$L_{loc}^{\infty}(\RR_+\times \RR,\, {\Cal H})$ (recall that ${\Cal H}
={\Cal H}_1\cup {\Cal H}_3$) which satisfy
$$
\hbox{system (1.1) in the sense of distributions,}
\leqno{(1.5a)}
$$
$$
\hbox{ the entropy inequality (0.2), (0.3) in the sense of distributions}
\leqno{(1.5b)}
$$
and
$$
\hbox{ an initial condition } u_0 \hbox{ at } t=0 \hbox{ in the } L_{loc}^1
\hbox{ sense.}
\leqno{(1.5c)}
$$
Here $u_0$ is a given function in $L_{loc}^\infty(\RR,{\Cal H})$ and,
for future reference, we rewrite the entropy inequality in the form
$$
\partial _t U(u) +\partial_x F(u)\le 0
\leqno{(1.6a)}
$$
with
$$
U(u) = W(w) +{v^2\over 2},\quad
F(u) = -\sigma (w)v \quad\hbox{ and }\quad
W(w) =\int_0^w \sigma(y)\, dy.
\leqno{(1.6b)}
$$
We recall that solutions in the sense (1.5) are {\it unique} (at least
for Riemann data) in the standard situation of a (genuinely nonlinear or 
linearly degenerate) {\it increasing} 
strain-stress function $\sigma$.  This is no longer true in the case
of the mixed system under consideration here: see for instance James [23]. 
We also point out that the entropy function $U$ is not a convex function.

To complete the formulation (1.5), we follow Abeyaratne and Knowles in [2]. 
Let us give first some
motivation for their suggestion. Suppose $u^\epsilon = (v^\epsilon,
w^\epsilon)$ is the solution of a regularized version of system (1.1)
obtained by adding high-order terms, depending on a (small) parameter
$\epsilon$, in the right-hand side of the equations (e.g. use
the viscosity-capillarity terms as was done by Slemrod [41]). As was
pointed out by Lax for general systems of conservation laws, the limit $u =\lim u^\epsilon$ -- 
if it exists (and if the convergence holds in a suitable topology) -- must be a solution to (1.1) in the sense (1.5);
in particular the entropy inequality (1.6) must hold. Since (1.5)
is incomplete, it seems natural to ``keep more information'' about the limiting function $u$
from its regularization $u^\epsilon$. Specifically  Abeyaratne and
Knowles' suggestion is equivalent to replacing (1.6) with the stronger requirement that
$$
\partial_t U(u) +\partial_x F(u) =\mu,
\leqno{(1.7)}
$$
where $\mu$ is a given non-positive measure that clearly must satisfy certain
restrictions. Note that in principle $\mu$ could be determined by the formula 
$$
\mu =\hbox{ weak-star }\lim_{\epsilon \to 0} (\partial _t U(u^\epsilon)
+\partial_x F(u^\epsilon)\, )
$$
(at least when $u^\epsilon$ has uniformly bounded total variation in $(t,x)$).
This formula may not give a very explicit expression for $\mu$.  Hopefully, it
turns out that (1.7) is needed (to achieve uniqueness) {\it only} for one
kind of discontinuity: the subsonic phase boundaries. Moreover, in that
case, we can allow a large range of measures $\mu$. Here, we call
subsonic (respectively supersonic) those phase boundaries that travel
with speed less (resp.~greater) than the contact discontinuities
in phase ${\Cal H}_3$.

The precise formulation of condition (1.7) given below requires that $u$
is a bounded function of bounded variation. When $u$ has bounded
variation, we call {\it entropy dissipation} the value of the measure
$\partial_t U(u) +\partial_x F(u)$ along a curve of (contact or phase)
discontinuity of $u$. According to [2], the {\it kinetic relation} yields
this entropy dissipation along any {\it subsonic phase boundary}, as an
explicit function, say $\phi(V)$, of the speed $V$ of propagation 
of this discontinuity. In applications, the actual kinetic
relation, that is the function $\phi$, must be determined from the properties 
of the specific material under consideration.
This kind of constitutive model is already in extensive
use in the quasi-static setting for problems of phase transition in solids.
We refer the reader to [1] as well as Truskinovsky [44] and the
references cited there.  The speed $V$ can also be interpreted as an
internal variable and the kinetic relation indeed determines the evolution
of this internal parameter.

\

\noindent{\bf Remark 1.1.} 1) That subsonic and supersonic phase boundaries must be treated
in a different way is clear, for instance when solving Riemann problems. 
A wave structure with a supersonic
phase boundary contains {\it two} waves, while one with a subsonic
boundary is composed of {\it three} waves. This latter case 
suffers, without a kinetic relation, from a strong lack of uniqueness. 
Cf.~James [23] and Section~2. 

\noindent 2) The approach considered here has some similarity to the theory
of nonlinear hyperbolic systems in non-conservative form; cf. Dal~Maso-LeFloch-Murat
[8] and LeFloch-Liu [30]. Namely, as is the case for systems (1.1), the weak solutions to these systems 
are not uniquely determined by the partial differential equations and
an entropy inequality, but an additional constitutive relation must be
added to ensure uniqueness. This fact was first pointed out by LeFloch; cf.~[27] to [29].

\noindent 3) Conservation laws with measure source-term  like (1.7) have been useful
in various contexts, cf.~Di~Perna [9], Di~Perna-Majda [11], Hou-LeFloch [21].

\vskip.2cm

Let us introduce some notations and recall some facts about functions
of bounded variation, that can be found in Volpert [45] and Federer
[14]. Let $\Omega$ be an open subset of $\RR^m$. A function $u\colon\,
\Omega\to\RR^p$ belongs to the space $BV(\Omega,\RR^p)$ (respectively
$BV_{loc}(\Omega,\RR^p)$) if $u\in L^1(\Omega,\RR^p)$ (resp.
$L^1_{loc}(\Omega,\RR^p)$) and the distributional derivatives
$\displaystyle{\partial u\over\partial y_j}$ for $1\le j\le m$
are bounded (resp. locally bounded) Borel measures on $\Omega$.
In what follows, we will always consider functions in $L^\infty(\Omega,
\RR^p)\cap BV(\Omega,\RR^p)$ or $L_{loc}^\infty(\Omega,\RR^p)\cap
BV_{loc}(\Omega,\RR^p)$, often called for short BV functions or 
$BV_{loc}$ functions. For each $BV_{loc}$ function $u$, we have
the following decomposition
$$
\Omega = C(u) \cup S(u)\cup E(u),
$$
where

$C(u)$ is the set of all points of approximate continuity for $u$,

$S(u)$ is the set of all points of approximate jump for $u$

\noindent and

$E(u)$ is the set of exceptional points with the property $H_{m-1}(E(u)\, ) = 0$. 

\noindent
Here $H_{m-1}$ is the 
$(m-1)$-dimensional Hausdorff measure on $\RR^m$. For each point $y$
in $S(u)$, there exists a unit normal $\nu\in\RR^m$ and approximate
left and right limits for $u$ that we denote by $u_\pm(y)$. The set $S(u)$ 
consists of the union of a countable number of rectifiable curves.

We denote the norm of $u$  by 
$\Vert u\Vert _{BV(\Omega,\RR^p)} =\Vert u\Vert_{L^1( \Omega,\RR^p)} 
+ |Du|(\Omega)$, where $Du$ is the measure
$(\displaystyle{\partial u\over\partial y_1}$,
$\displaystyle{\partial u\over\partial y_2}$, $\ldots$,
$\displaystyle{\partial u\over\partial y_m}$).
When $u = u(t,x)\in L_{loc}^\infty(\RR_+\times \RR, {\Cal H})\cap
BV_{loc}(\RR_+\times \RR, {\Cal H})$, we use the notation:
$$
\nu(t,x) = (\nu_t(t,x), \, \nu_x(t,x)\, )
\quad\hbox{ and }\quad 
V(t,x) = - {\nu_t(t,x)\over \nu_x(t,x)}
$$
valid for all $(t,x)\in S(u)$. The ratio $V(t,x)$ represents the speed
of propagation of the discontinuity in $u$ at the point $(t,x)$.
Note that system (1.1) has the property of propagation with finite
velocity (in regions ${\Cal H}_1$ and ${\Cal H}_3$). So $\nu_x(t,x)$
will never vanish, and for definiteness we always choose
$\nu_x(t,x) > 0$. In the following, we shall always have: $u(t) \in BV$ for all times $t$.

Let $\phi\colon\, ]-c_3,c_3[ \to \RR$ be a function, called below {\it entropy
dissipation function}, satisfying the following properties:
$$
\phi \hbox{ belongs to } {\Cal C}^2 (]-c_3,0[\cup ] 0,c_3]) \,\,
\hbox{ and }\,\, \phi(0\pm) \hbox{ and }\,\, \phi'(0\pm) \hbox{ exist,}
\leqno{(1.8a)}
$$
$$
\lim_{V\to c_3^-} \phi =\overline \psi (c_3) 
\quad\hbox{ and }\quad
\phi'''(c_3-) \hbox{ exists, } 
\leqno{(1.8b)}
$$
$$
\lim_{V\to-c_3^+} \phi =-\infty,
\leqno{(1.8c)}
$$
$$
\phi \hbox{ is increasing on } ]-c_3,c_3]
\leqno{(1.8d)}
$$
and
$$
\cases
\underline{\psi}(V) \le \phi(V) \le 0\hfil & \text{ for }\quad V\in \,]-c_3,0], \cr
  \cr
0 \le \phi(V) \le \overline \psi(V)\hfil & \text{ for }\quad V\in [0,c_3].\cr
\endcases
\leqno{(1.8e)}
$$
In (1.8b) and (1.8e), the {\it minimal} and {\it maximal entropy dissipation
functions} $\underline{\psi}\colon\, ] -c_3,0]\to\RR_-$ and
$\overline\psi\colon\, [0,c_3] \to \RR_+$ are defined by
$$
\underline{\psi}(V) = {(k_1-k_3)\over 2} w_M(w_m - {k_1-V^2\over k_3-V^2}
w_M) \quad\hbox{ for }\quad V \in ]-c_3,0]
\leqno{(1.9a)}
$$
and
$$
\overline{\psi}(V) = {(k_1-k_3)\over 2} w_m(w_M - {k_3-V^2\over k_1-V^2}
w_m) \quad\hbox{ for }\quad V \in [0,c_3].
\leqno{(1.9b)}
$$

\vskip.2cm

\noindent{\bf Remark 1.2.} 1)~Inequalities (1.8e) give the range of values taken 
by the entropy dissipation rate ${\Cal E}(u)$ (see below) when varying the left 
and right values at a discontinuity satisfying the Rankine-Hugoniot relations 
and the entropy condition.

\noindent 2)~In [2], instead of (1.8e), Abeyaratne and Knowles assume the 
(more restrictive) condition:
$$
\cases
\underline{\psi}(0)\le \phi(V)\le 0 \hfil & \text{ for } \quad V\in \,]-c_3,0], \cr
  \cr
\phantom{x}0\le \phi(V)\le \overline\psi(0) \hfil & \text{ for }\quad V\in [0,c_3].\cr
\endcases
\leqno{(1.8e)'}
$$

\noindent 3)~Assumptions (1.8) made in this paper are indeed satisfied in the 
examples considered by [3] and [4]. For instance, they are fulfilled by
 the {\it maximally dissipative function\/} $\phi_{\text{max}}$ defined by:
$$
\phi_{\text{max}}(V) =\underline{\psi}(V) \quad\hbox{ for }\quad
V \in \,]-c_3,0], \quad \overline{\psi}(V)\quad\hbox{ for }\quad
V\in [0,c_3].
$$
\vskip.2cm

%


\

We next define the {\it entropy dissipation rate} ${\Cal E}(u)$
associated with any function 
$$
u\in L_{loc}^\infty (\RR_+\times \RR,
{\Cal H}) \cap BV_{loc}(\RR_+\times \RR,{\Cal H})
$$
 by the following
formula
$$
{\Cal E} (u) = - (U(u_+) - U(u_-)\, ) -{\nu_x\over\nu_t} (F(u_+)
- F(u_-)\, )
\leqno{(1.10)}
$$
which defines ${\Cal E}(u)(t,x)$ at $H_1$--almost every point $(t,x)$,
where $\nu_t(t,x)\ne 0$ (i.e. $V(t,x)\ne 0$).
${\Cal E}(u)$ is the product of $-{1\over\nu_t}$ by the jump of the
measure $\partial_tU(u) +\partial_xF(u)$ along the curve of approximate
jump of $u$. Formula (1.10) makes sense only if $\nu_t(t,x)\ne 0$.
However, it is a simple observation that if $u$ is assumed to be a
{\it weak solution} to system (1.1), then the above jump (i.e. the entropy
dissipation) vanishes at the points where $\nu_t$ vanishes.
This fact allows us to define ${\Cal E}(u)(t,x)$ $H_1$--almost everywhere, as shown by the following lemma.

\

\noindent{\bf Lemma 1.1.} {\it If} $u\in L_{loc}^\infty(\RR_+ \times \RR, {\Cal H})
\cap BV_{loc}(\RR_+\times \RR, {\Cal H})$ {\it is a weak solution to} (1.1),
{\it then one has}
$$
{\Cal E}(u) =- \int_{w_-} ^{w_+} \bigg\{\sigma (y)-{1\over 2} (\sigma(w_+)
+ \sigma(w_-)\, )\bigg\} \, dy,
\leqno{(1.10)'}
$$
{\it at} $H_1${\it --almost every} $(t,x)$ {\it such that} $\nu_t(t,x) \ne 0$.

\vskip.2cm

From now on, we use (1.10)$'$ to define ${\Cal E}(u)(t,x)$.

\

\noindent{\bf Proof.} At a point of approximate discontinuity $(t,x)$ of the
solution $u$, the following Rankine-Hugoniot relations hold:
$$
\cases
\nu_t(w_+ - w_-) -\nu_x(v_+-v_-) = 0, \cr
   \cr
\nu_t (v_+ - v_-) - \nu_x(\sigma(w_+)-\sigma(w_-)\, ) = 0. \cr
\endcases
$$
These relations used in (1.10) yield:
$$
\eqalign{
-{\Cal E} (u) &= \int_{w_-}^{w_+} \sigma(y)\, dy + {1\over 2} (v_+^2
- v_-^2) - {\nu_x\over\nu_t} (\sigma(w_+) v_+ - \sigma(w_-)v_-) \cr
&=\int_{w_-}^{w_+} \sigma(y)\, dy + {1\over 2} (v_+ + v_-) {\nu_x
\over \nu_t} (\sigma(w_+ ) -\sigma(w_-)\, ) -{\nu_x\over \nu_t}
(\sigma(w_+)v_+ - \sigma(w_-)v_-). \cr
}
$$
We thus get
$$
\eqalign{
-{\Cal E}(u) &= \int_{w_-}^{w_+} \sigma(y)\, dy + {\nu_x\over 2\nu_t}
\big\{ v_+\sigma (w_+) + v_- \sigma(w_+) - v_+ \sigma(w_-) - v_-
\sigma(w_-) \cr
&\quad - 2 \sigma(w_+)v_+ + 2\sigma (w_-)v_-\big\}, \cr
}
$$
so that
$$
-{\Cal E} (u) =\int_{w_-}^{w_+} \sigma(y)\, dy -{1\over 2} \, {\nu_x
\over \nu_t} (v_+-v_-) (\sigma(w_+)+\sigma(w_-)\, ),
$$
which, in view of the Rankine-Hugoniot relations above, gives the desired result (1.10)$'$. ~~$\Box$

\vskip.2cm

Let us denote by $B_{sub}(u)$ the set of all points of approximate
discontinuity in a weak solution $u$ that correspond to a {\it subsonic}
phase boundary. This means:
$$
\eqalign{
B_{sub}(u) &= \bigg\{ (t,x) \in S(u) | \hbox{ either : } u_- (t,x)
\in {\Cal H}_1, \, u_+ (t,x)\in {\Cal H}_3 \hbox{ and } |V| \le c_3,\cr
&\qquad\qquad \hbox{ or : } u_-(t,x)\in {\Cal H}_3, \, 
u_+ (t,x)\in {\Cal H}_1 \hbox{ and } |V| \le c_3\bigg\}.\cr
}
$$
In view of (1.6), the Borel measure $\partial_t U(u) + \partial_x F(u)$
is globally non-positive. The {\it kinetic relation\/} now specifies the value
itself (and not only the sign) of this measure along any subsonic phase
boundary. In other words, for $H_1$--almost all $(t,x)\in {\Cal B}_{sub}(u)$,
one must have
$$
{\Cal E}(u)(t,x) = \cases
-\phi(V(t,x)\, ) \hfil & \text{ if }\quad u_-(t,x)\in {\Cal H}_1, \cr
  \cr
\phantom{-}\phi(-V(t,x)\, ) \hfil & \text{ if } \quad 
u_-(t,x) \in {\Cal H}_3. \cr
\endcases
\leqno{(1.11)}
$$
\vskip.2cm

\noindent{\bf Remark 1.3.} As a matter of fact, the traveling waves obtained 
through the viscosity-capillarity regularization to system (1.1)
converge to weak solutions of (1.1) that satisfy the kinetic relation
(1.1) with a specific choice of function $\phi$. This function can be 
determined explicitly and depends only on the viscosity and capillarity
coefficients introduced in the regularization (cf. [3]).

\vskip.2cm

Finally, we have to formulate the {\it initiation criterion}, which
together with the above kinetic relation will allow us to rule out
all non-physical solutions to our problem. Let $]a,b[$ be a space interval
 in which we are going to set the
problem, with $a < b$ and possibly $a=-\infty$ and/or $b=+\infty$.
The initiation criterion will reflect the following facts:
$$
\left\{
\eqalign{
&\hbox{ no new phase } \hbox{occurs from any point } x \hbox{ in } ]a,b[\cr
&\hbox{ except if no solution exists without creation of a 
new phase,}\cr
}
\right.
\leqno{(1.12)}
$$
$$
\left\{
\eqalign{
&\hbox{a new phase state may occur at the boundary point } x=a, \cr
&\hbox{ even if a solution with no new phase exists;}\cr
 &\hbox{ a criterion is required to make the choice}\cr
}
\right.
\leqno{(1.13)}
$$
and
$$
\left\{
\eqalign{
&\hbox{a new phase state may occur at } x=b, \cr
&\hbox{ even if a solution with no new phase exists;} \cr
&\hbox{ a criterion is required to make the choice.}\cr
}
\right.
\leqno{(1.14)}
$$

From the mathematical point of view, condition (1.12) is essential: it
ensures that spontaneous initiation of a new phase inside $]a,b[$ cannot
occur from two nearby initial states in the same phase (cf.~Section~2).
This does not exclude the possibility (and it really happens) that an
initial discontinuity with large jump gives rise to, for instance,
a phase~1 state although the states on both sides of the initial
discontinuity are in phase~3. However, by condition (1.12), a single
{\it constant state\/} is always a (trivial) {\it admissible solution\/}. 
(This property was not satisfied in the construction of [2].) This is also essential to get 
the $L^1$ continuous dependence property for Riemann solutions, 
proved below in Section~2.

Conditions (1.13) and (1.14) follow the
quasi-static theory [1]. They allow ``spontaneous nucleation'' of a new phase 
only at the end points of $[a,b]$. Note that, more generally, we could as well
 allow nucleation at some arbitrary given points of $[a,b]$. Our actual restriction 
is that the points of spontaneous nucleation are known a priori and follow a 
selection criterion of the form specified below. However, while this formulation is fully 
satisfactory from the mathematical point of view, it does not 
reproduce what is really observed in practical experiments with
elastic bars. Namely, in experiments, when pulling out an elastic
bar uniformly in phase~1, initiation of phase~3 regions
in the bar occurs successively and (apparently) randomly at various
places in the bar. Physicists assert that initiation occurs at microscopic 
inhomogeneities of the material. A complete treatement of the initiation 
mecanism is beyond the scope of this paper and would probably require a 
statistical description. (As a matter of fact, this might be included 
in the random choice scheme, studied below, quite easily.)

It remains to provide an analytic version of the conditions (1.12) to (1.14).
For convenience, we use here an averaged strain in our formulation.
(In [1] and [2], the stress and the entropy dissipation
rate, respectively, are used instead.) Given any function $u= (v,w)$ in $L_{loc}^\infty
(\RR_+\times \RR, {\Cal H}) \cap BV_{loc}(\RR_+\times \RR, {\Cal H})$,
we set
$$
h_u = {c(w_-) w_- + c(w_+) w_+ + v_+ - v_- \over c(w_-)+c(w_+)}
\leqno{(1.15)}
$$
which defines $h_u(t,x)$ for $H_1$ -- almost every $(t,x)$ in
$\RR_+ \times\RR$. We note that $h_u(t,x) = w(t,x)$ when $(t,x)$ is
a point of approximate continuity of $u$. So $h_u(t,x)$ represents
an {\it averaged strain} at the point $(t,x)$ and determine the dynamics at this point. 
For instance, if $u_-$ and $u_+$ are in the same phase, $h_u$ is the intermediate value
 between the 1-wave and the 2-wave in the solution of the Riemann problem with initial 
data $u_-$ and $u_+$ (cf. Section 2).

At any interior point $x\in ]a,b[$ and for each time $t \ge 0$, the initiation criterion, 
by definition,
is
$$
\left\{
\eqalign{
&\hbox{if } u_-(t,x) \hbox{ and } u_+(t,x) \hbox{ belong to } {\Cal H}_1 ~
(\hbox{respectively } {\Cal H}_3), \hbox{ then:} \cr
&\quad h_u (t,x) \in {\Cal H}_1 ~(\hbox{resp. } {\Cal H}_3) 
\hbox{ if and only if there exists } \epsilon > 0 \hbox{ such that} \cr
&u(s,y)\in {\Cal H}_1 ~ (\hbox{resp. } {\Cal H}_3) 
\hbox{ for } (s,y)\in[t,\, t+\epsilon [\times ]x-\epsilon, x+\epsilon [. \cr
}
\right.
\leqno{(1.16)}
$$
According to (1.12), condition (1.16) ensures that, locally in time,
the solution remains in the same phase whenever this is possible.
Cf.~Section~2.

We are now concerned with the boundary points $x=a$ and $x=b$. We assume that $u(t)$ 
is defined for all times and has bounded variation in $x$. (This will be the regularity
 of the solutions found in Section 4.) The material is assumed to be fixed at the end
 points, i.e. when $a\ne -\infty$ and/or $b\ne +\infty$, we have
$$
v_+ (t,a) = 0 \quad\hbox{ for } L^1\hbox{-- almost every } t > 0,
\leqno{(1.17a)}
$$
and
$$
v_- (t,b) = 0 \quad\hbox{ for } L^1\hbox{-- almost every } t > 0,
\leqno{(1.17b)}
$$
where $L^1$ denotes the one-dimensional Lebesgue measure. Since $v$ has bounded variation, 
it admits a $L^1$ trace at $x=a$ and $x=b$. Let $w_M^{cr}$
and $w_m^{cr}$ be two constants, called {\it critical values for the
initiation}, that must satisfy the inequalities
$$
{\sigma_0\over k_1} \le w_M^{cr} \le w_M \quad\hbox{ and }\quad
w_m \le w_m^{cr} \le {\sigma_0\over k_3},
\leqno{(1.18)}
$$
where $\sigma_0$ is the so-called {\it Maxwell stress} given by
$$
\sigma_0 =\sqrt{\sigma_m\sigma_M} = c_1 c_3\sqrt{w_mw_M}.
\leqno{(1.19)}
$$
Note that, as pointed out to us by Abeyaratne, the critical values for intitiation should
 in principle depend on the speed of propagation of the phase discontinuity. 
At the point $x=a$ (when $a\ne -\infty$), we impose for {\it all} times
$t\ge 0$ the following two conditions:
$$
\left\{
\eqalign{
&\hbox{if } u_+(t,a) \hbox{ belongs to } {\Cal H}_3, \, \hbox{ then:}\cr
&h_u(t,a)\ge w_m^{cr} \hbox{ if and only if there exists }
\epsilon > 0 \hbox{ such that } \cr
&u_+(s,a)\in {\Cal H}_3 \hbox{ for } s\in [t,t+\epsilon[, \cr
}
\right.
\leqno{(1.20)_i}
$$
and
$$
\left\{
\eqalign{
&\hbox{if } u_+(t,a) \hbox{ belongs to } {\Cal H}_1 \, \, \hbox{ then:}\cr
&h_u(t,a)\in {\Cal H}_1 \hbox{ if and only if there exists } \epsilon > 0
\hbox{ such that} \cr
&u_+(s,t) \in {\Cal H}_1 \hbox{ for all } s\in [t,t+\epsilon
[. \cr
}
\right.
\leqno{(1.20)_{ii}}
$$
In order to satisfy the boundary condition (1.17a), the term $h_u(a,t)$ 
in (1.20)$_i$ is defined by formula (1.15) with here
$$
v_-(t,a) = -v_+ (t,a) \quad\hbox{ and }\quad w_-(t,a) = w_+(t,a).
\leqno{(1.21)}
$$

Similarly, at the point $x=b$ (when $b\ne +\infty$), we impose for {\it all}
times $t>0$ that
$$
\left\{
\eqalign{
&\hbox{if } u_-(t,b) \hbox{ belongs to } {\Cal H}_1, \, \hbox{then:}\cr
&h_u(t,b) < w_M^{cr} \hbox{ if and only if there exists } \epsilon > 0
\hbox{ such that }\cr
&u_-(s,b)\in {\Cal H}_1\, \hbox{ for all }
s\in [t,t+\epsilon[, \cr
}
\right.
\leqno{(1.22)_i}
$$
and
$$
\left\{
\eqalign{
&\hbox{ if } u_-(t,b) \hbox{ belongs to } {\Cal H}_3, \, \hbox{then: }\cr
&h_u(t,b)\in {\Cal H}_3 \hbox{ if and only if there exists } \epsilon > 0
\hbox{ such that } \cr
&u_-(s,b)\in {\Cal H}_3 \hbox{ for all } s \in [t,t+\epsilon[. \cr
}
\right.
\leqno{(1.22)_{ii}}
$$
As previously, we set here:
$$
v_+(t,b) = - v_-(t,b) \quad\hbox{ and }\quad w_+ (t,b) = w_-(t,b).
\leqno{(1.23)}
$$

We call an {\it admissible weak solution} to system (1.1) a function
$u=(v,w)$ which satisfies the
{\it conservation laws } (1.1), the {\it entropy inequality } (1.6),
the {\it kinetic relation } (1.11), the {\it boundary condition } (1.17)
(if instead of $\RR$ an interval $]a,b[$ is considered) and the 
{\it initiation criterion } (1.16), (1.20) and (1.22).

In this paper, we will prove the existence of such an admissible weak
solution for two kinds of Cauchy data: the Riemann problem (in the whole space
and in a half space) and a perturbation of a single propagating phase
boundary. These results provide a strong justification for our formulation
here. It would be interesting to address the general question of existence
and uniqueness for system (1.1) in the setting introduced in this section.

\

\noindent{\bf Remark 1.4.} 1)~A phase boundary necessarily is a wave with a
{\it large} strength (at least $|w_M-w_m|$). This implies that,
in BV solutions, phase boundaries cannot accumulate in a bounded region
of the $(t,x)$--plane. Phase boundaries are thus isolated, and this
justifies the formulations (1.16), (1.20) and (1.22).

\noindent 2)~Sections~3 and 4 will provide an existence result for {\it small } BV
perturbations of phase boundaries. We believe our result to be true
for any {\it finite } number of phase boundaries. However, for {\it arbitrary large }
data, the appearance of an infinite number of phase boundaries is not
excluded a~priori.  In such a case, the solution
would not have bounded total variation. A challenging issue is to extend the present 
formulation to the framework of $L_\infty$ solutions.

\noindent 3)~The formulation of this section can be extended to the case
that the stress-function is not a piecewise affine function
but an arbitrary piecewise monotone function.

\noindent 4)~Shearer's solution [38] corresponds to the choice $w_M^{cr}=w_M$, 
$w_m^{cr}=w_m$ and $\phi=\phi_{\text{max}}$ (see Remark 1.2 for the definition of $\phi_{\text{max}}$).

\vskip.2cm

\vfill
\eject
\noindent{\bf 2. RIEMANN PROBLEM AND CONTINUOUS DEPENDENCE}

\vskip.2cm

This section gives an explicit description of the admissible weak solution of
the problem formulated in Section~1, in two cases: the Riemann problem
in the whole space and the Riemann problem in a half space. Our main result in
 this section is the $L^1$ continuous
dependence property of the solution of these problems. Note that the formulation 
of Section~1 and the assumptions made
there are essential for this property to hold.

We consider the following two problems:

1) {\it the Riemann problem in the whole space} $]a,b[ = ]-\infty,+\infty[$
which corresponds to an initial data of the form
$$
u_0(x) = u_L \quad\hbox{ for }\quad x < 0 \quad\hbox{ and }\quad
u_R \quad\hbox{ for }\quad x> 0.
\leqno{(2.1)}
$$
Here $u_L = (v_L,w_L) \in {\Cal H}$ and $u_R = (v_R,w_R) \in {\Cal H}$
are two given constant states.

2) {\it the Riemann problem in the half space} $]a,b[ = ]0,+\infty[$ which
corresponds to the initial data
$$
u_0(x) = u_0 \quad\hbox{ for all }\quad x > 0
\leqno{(2.2)}
$$
where $u_0 = (v_0,w_0 ) \in {\Cal H}$ is a constant state. 

\noindent We shall describe successively the admissible solutions to problems 1) and 2), by following
closely the work by Abeyaratne-Knowles. However our
construction is slightly different from the one in [2], due to
our formulation.  We will not address the question of
uniqueness of the solution here, since it is an easy matter from the
results in [2] (which yield for their construction uniqueness 
in the class of solutions composed of simple waves).

To begin with, we deal with problem 1) and distinguish between several cases:

case 1)-a$\colon \,\quad u_L \in {\Cal H}_1 \hbox{ and } u_R\in {\Cal H}_3$,

case 1)-b$\colon \,\quad u_L \in {\Cal H}_1 \hbox{ and } u_R\in {\Cal H}_1$,

case 1)-c$\colon \,\quad u_L \quad {\Cal H}_3 \hbox{ and } u_R\in {\Cal H}_3$,

case 1)-d$\colon \,\quad u_L \in {\Cal H}_3 \hbox{ and } u_R\in {\Cal H}_1$.

\noindent Cases  1)-c and 1)-d are very similar to cases
1)-b and 1)-a, respectively. (Use the transformation $x\to- x$ and
the fact that the equations (1.1) and more generally all the requirements 
in the formulation of Section~1 are invariant under this transformation.) 
So we omit them and focus on the two first cases.

\vskip.2cm

\noindent{\bf Case 1)-a:} Suppose that $u_L \in {\Cal H}_1$ and
$u_R\in {\Cal H}_3$.

\vskip.2cm

We must construct a solution to (1.1), (2.1) which is admissible in the
sense of Section~1. The solution necessarily contains a phase boundary
(and only one as was checked in [2]) with phase~1 states at the left and
phase~3 states at the right.  Two different wave structures are possible,
depending on whether the phase boundary is {\it subsonic} or {\it supersonic}.
Let $V$ be the speed of the phase boundary and set
$$
h_{L R} = h_u(0,x) = {1\over c_1 + c_3} (c_1 w_L + c_3 w_R + v_R-v_L).
\leqno{(2.3)}
$$
We distinguish between two cases depending on the sign of $h_{LR}$.

\

\noindent{\bf Case 1)-a1:} Suppose moreover that $h_{LR} > 0$.

In that case, we seek the solution $u$ in the following form
$$
u(t,x) = \cases
u_L \hfil & \text{ for } \quad x < -c_1 t, \cr
u_- \hfil & \text{ for } \quad -c_1t < x < Vt, \cr
u_+ \hfil & \text{ for } \quad vt < x < c_3t\cr
u_R \hfil & \text{ for } \quad x>c_3t, \cr
\endcases 
\leqno{(2.4)}
$$
where the constants $u_- = (v_-,w_-)$ and $u_+ = (v_+,w_+)$ belong to
${\Cal H}_1$ and ${\Cal H}_3$ respectively.
The
solution contains a contact discontinuity of speed $-c_1$, the phase
boundary with subsonic speed $|V|<c_3$ and a contact discontinuity with
speed $c_3$. We are going to prove that indeed such a solution exists
by determining explicitly the values of the constants $u_-$, $u_+$
and $V$.

First of all, $u$ given by (2.4) must be a weak solution to (1.1),
so must satisfy the following four Rankine-Hugoniot relations
$$
\left\{
\eqalign{
c_1 (w_- - w_L) - (v_- - v_L) = 0, \qquad
&V(w_+ - w_-) + v_+ - v_- = 0, \cr
V(v_+ - v_-) + c_3^2 w_+ - c_1^2w_- = 0, \qquad
&c_3(w_+ - w_R) + v_+ - v_R = 0. \cr
}
\right.
\leqno{(2.5)}
$$
If $V$ is chosen as a parameter, (2.5) yields explicit expressions
for $v_-$, $v_+$, $w_-$, $w_+$ as functions of $V$:
$$
\left\{
\eqalign{
v_- = v_L - c_1 w_L + {c_3+V\over c_1+V}\, c_1 h_{LR}, \qquad
&w_- = {c_3+V \over c_1 + V} \, h_{LR},\cr
v_+ = v_R + c_3 w_R - {c_1-V \over c_3-V} \, c_3 h_{LR}, \qquad
&w_+ = {c_1-V\over c_3-V} \, h_{LR}, \cr
}
\right.
\leqno{(2.6)}
$$
Formulas (2.6) define a {\it one-parameter family of solutions} to
problem (1.1), (2.1). Note that $w_-$ and $w_+$ are always non-negative.

Next we take into account the kinetic relation that states (cf. (1.10)
and (1.11)):
$$
U(u_+) - U(u_-) + {1\over V} \big( F(u_+) - F(u_-)\, \big) =\phi(V),
$$
or using the more general form (1.10)$'$:
$$
\int_{w_-} ^{w_+} \big\{\sigma (y) -{1\over 2} \big( \sigma(w_+)
+\sigma(w_-)\, \big)\, \big\} \, dy = \phi(V).
$$
Using the expression (1.2) for the function $\sigma$, this becomes
$$
{1\over 2} (k_1 - k_3) \, (w_M w_m - w_+w_-) = \phi(V).
\leqno{(2.7)}
$$
If we use in (2.7) the expressions for $w_+$ and $w_-$ given by (2.6),
it follows that
$$
\theta(V) =\phi(V),  \quad\hbox{ where }\quad \theta(V) = {(k_1-k_3)
\over 2} \bigg\{ w_M w_m - {(c_3+V)(c_1-V)\over (c_1+V)(c_3-V)} \, h^2_{LR}
\bigg\}.
\leqno{(2.8)}
$$
Note that the function $\theta$ {\it depends only on the averaged strain}
$h_{LR}$. In view of our set of assumptions (1.8) and 
$$
\theta ' < 0, \quad \theta(-c_3)=\overline \psi(c_3) \quad\hbox{ and }\quad
\lim_{V\to c_3} \theta(V) = -\infty,
$$
one easily checks that
equation (2.8) admits a unique root $V$. 
Moreover, if this specific value of $V$ is used in (2.6)
to get $w_-$, $w_+$, $v_-$ and $v_+$, then our construction is consistent 
in the sense that one has
$$
w_-\in {\Cal H}_1 \quad\hbox{ and }\quad w_+ \in {\Cal H}_3.
\leqno{(2.9)}
$$
We now prove (2.9) by using the assumption (1.8e) on the function $\phi$ (a stronger
assumption was made in [2] to derive (2.9)).

Let us for instance check that $w_-\in {\Cal H}_1$, in other words:
$w_- \le w_M$. In view of (1.8e) and (2.8), one has $\underline{\psi}(V)
< \theta(V)$, so that by (1.9a) and (2.8) again
$$
w_M w_m - {k_1-V^2\over k_3-V^2} \, w_M^2 < w_M w_m - {(c_3+V)(c_1-V)
\over (c_1+V)(c_3-V)} \, h^2_{LR}.
$$
Since $|V| < c_3$, we obtain
$$
w_M^2 > h^2_{LR} \, {(c_1+V)^2\over (c_3+V)^2}.
$$
But $w_M > 0$ and $h_{LR} > 0$ by assumption, so
$$
w_M > h_{LR} \, {c_1+V\over c_3+V} = w_-
$$
in view of the expressions (2.6). The proof of (2.9) is complete. 

Finally, we note that the entropy inequality (1.6) is trivially
satisfied along the contact waves, while it is a consequence of the 
kinetic relation (1.10), (1.11) along the phase boundary. Thus, in the 
present case, (1.6) yields no additional constraint.

The above construction yields the admissible weak solution of the
problem.  Based on the explicit expressions (2.6) and the implicit
equation (2.8), it is elementary to prove the following regularity
result of the Riemann solution.

\ 

\noindent{\bf Lemma 2.1.} {\it Consider the Riemann problem} (1.1), (2.1)
{\it in case} 1)-a1, {\it i.e. with} $u_L \in {\Cal H}_1$, $u_R\in {\Cal H}_3$
{\it and} $h_{LR}> 0$. {\it Then the admissible weak solution to this
problem is given by formulas} (2.4), (2.6) {\it and} (2.8). {\it One can
consider the states} $u_-$ {\it and} $u_+$ {\it and the speed} $V$ {\it in} 
(2.4) {\it as functions of the initial states} $u_L$ {\it and} $u_R$, 
{\it or more precisely}
$$
v_- = v_-(u_L, h_{LR}), \qquad w_- = w_- (h_{LR}) 
$$
{\it and}
$$
v_+ = v_+(u_R, h_{LR}), \quad w_+ = w_+(h_{LR}) \quad\hbox{\it and }\quad
V=V(h_{LR}).
$$
{\it The functions $u_-$, $u_+$ and $V$ are Lipschitz continuous 
in the range of values $\{ u_L\in {\Cal H}_1, \, u_R\in {\Cal H}_3
\, | \, h_{LR} > 0\}$. They are of class  ${\Cal C}^1$ (with Lipschitz continuous 
derivatives) away from $V=0$. The behavior of $u_-$, $u_+$ and $V$ when $h_{LR}\to 0^+$
is given as follows:} 
$$
\leqalignno{
&\lim_{h_{LR}\to 0+} v_- = v_L - c_1 w_L, 
\phantom{xxxxxxxxxxx}
\lim_{h_{LR}\to 0+} w_- = 0, &(2.10a)\cr
&\lim_{h_{LR} \to 0+} v_+ = v_R + c_3\bigg\{ w_R - \sqrt{ {\phi'(c_3)
\over c_3} }\bigg\}, 
\lim_{h_{LR}\to 0+} w_+=\sqrt{ {\phi'(c_3)
\over c_3} }, &(2.10b)\cr
&\lim_{h_{LR}\to 0+} V= c_3,  
\phantom{lxxxxxxxxxxxxxxxxx}
\lim_{h_{LR}\to 0+} {dV\over dh_{LR}} = (c_3 - c_1) \sqrt{ {c_3
\over \phi '(c_3)} }  &(2.10c)\cr
& &\hbox{ and } \cr
&\lim_{h_{LR} \to 0+} {\partial v_-\over \partial h_{LR}} = {2c_1c_3
\over c_1 + c_3}, 
\phantom{xxxxxxxxxxx}\lim_{h_{LR}\to 0+} {dw_-\over dh_{LR}} =
{2c_3\over c_1+c_3}. &(2.10d)\cr
}
$$

\vskip.2cm

\noindent{\bf Remark 2.1.}  1)~Assumptions (1.8a) and (1.8d) imply that $\phi^{-1}$
exists and is a Lipschitz continuous function. Away from $V=0$,  $\phi^{-1}$ is of class 
  ${\Cal C}^2$, and so is the function $V(h_{LR})$
in view of (2.8). 

\noindent 2)~If the function $\phi \in {\Cal C}^2(]-c_3,c_3])$, then all the functions 
in Lemma 2.1 are globally of class ${\Cal C}^2$. (This is {\it not\/} the case of
 the maximally dissipative function quoted in Remark 1.2.)

\noindent 3)~In the special case that $\phi=\phi_{\text{max}}$, we find: $\sqrt{ {\phi'(c_3)
\over c_3} } = w_m$.

\vskip.2cm

\ 



\noindent{\bf Case 1)-a2:} Suppose now that $h_{LR}\le 0$.

In that case, the solution is composed of a contact discontinuity with
speed $-c_1$ and a phase boundary with {\it supersonic} speed $V>c_3$.
There is no $c_3$ contact wave. We use the notation
$$
u(t,x) = \cases
u_L \hfil & \text{ for }\quad x< -c_1t, \cr
u_- \hfil & \text{ for }\quad -c_1t < x < Vt, \cr
u_R \hfil & \text{ for }\quad x>Vt.\cr
\endcases
\leqno{(2.11)}
$$
The state $u_- = (v_-, w_-) \in {\Cal H}_1$ and the speed $V$
must satisfy the following jump relations:
$$
\eqalign{
&c_1(w_--w_L) - (v_- -v_L) = 0, \cr
&V(w_R-w_-) + v_R-v_- = 0, \qquad V(v_R-v_-)+c_3^2w_R-c_1^2w_-=0 \cr
}
$$
We thus get $v_-$ and $w_-$ explicitly as functions of $V$
$$
v_- = v_L - c_1 w_L +{c_3+V\over c_1+V} \, c_1 h_{LR}, \qquad
w_- = {c_3+V\over c_1+V} \, h_{LR},
\leqno{(2.12a)}
$$
the speed $V$ being given by an implicit algebraic equation
$$
V^2\big\{ -c_3 w_R + (c_1 +c_3) h_{LR}\big\} + (c_3^2 - c_1^2)w_RV
+ (c_1+c_3) c_1(c_3w_R - c_1h_{LR}) = 0.
\leqno{(2.12b)}
$$
Note that $w_-$ given by (2.12a) is always non-positive.

One can check [2] that (2.12b) has a unique solution which
belongs indeed to the (physically interesting) interval $[c_3,c_1[$
if and only if $h_{LR}$ satisfies the restriction
$$
h_\infty < h_{LR} \le 0 \quad\hbox{ with }\quad h_\infty = {1
\over c_1 +c_3} \big\{ c_3 w_R - c_1 
- (c_3^2 w_R^2 + (c_1^2 + c_3^2)w_R + c_1^2)^{1/2} \big\}.
\leqno{(2.13)}
$$
In other words, the Riemann problem can be solved when $h_{LR}\le 0$
if and only if $h_{LR} > h_\infty$. We emphasize that the kinetic
relation was not used here and $V$ is found to be supersonic; this is in
complete agreement with the fact that (1.11) is imposed 
only for subsonic phase boundaries.

\

\noindent{\bf Lemma 2.2.} {\it Consider the Riemann problem} (1.1), (2.1)
{\it in the case} 1)-a2 {\it that} $u_L\in {\Cal H}_1$, $u_R\in {\Cal H}_3$
{\it and} $h_{LR}\le 0$ ({\it with the restriction} (2.13)). {\it Then the
admissible weak solution to this problem is given by} (2.11), (2.12).
{\it One can consider the state} $u_-$ {\it and the speed} $V$ {\it as 
functions of the initial states} $u_L$ {\it and} $u_R$, {\it or more 
precisely:}
$$
v_- = v_-( u_L, h_{LR}, w_R), \qquad
w_- = w_-(h_{LR}, w_R) \quad\hbox{ and }\quad
V = V(h_{LR}, w_R).
$$
\noindent {\it Then the functions} $u_-$ {\it and} $V$ {\it are} ${\Cal C}^\infty$ 
{\it functions of their arguments and when} $h_{LR}\to 0-$, {\it they satisfy:}
$$
\leqalignno{
&\lim_{h_{LR}\to 0-} v_- = v_L - c_1 w_L, \phantom{xx}\qquad
\lim_{h_{LR}\to 0-} w_- = 0, &(2.14a)\cr
&\lim_{h_{LR}\to 0-} V = c_3, \phantom{xxxxxxxx} \qquad
\lim_{h_{LR}\to 0-} {\partial V\over\partial h_{LR}} = {(c_1+c_3)^2(c_1-c_3)
\over (c_1^2 + c_3^2) w_R}, &(2.14b)\cr
& & \hbox{and} \cr
&\lim_{h_{LR} \to 0-} {\partial v_-\over\partial h_{LR}} = {2c_1c_3
\over c_1+c_3},\phantom{xx}\qquad 
\lim_{h_{LR}\to 0-} {dw_-\over dh_{LR}} = {2c_3 \over c_1+c_3}. &(2.14c)\cr
}
$$

\vskip.2cm

\noindent{\bf Remark 2.2.} The limits found in (2.14a) and (2.14c) coincide
with those in (2.10a) and (2.10d) respectively. This implies that,
{\it except} at those points where $V=0$, the function $u_- = (v_-,w_-)$
is of class ${\Cal C}^1$ (with Lipschitz continuous derivatives) in the {\it whole} domain $\{ u_L \in {\Cal H}_1$,
$u_R\in {\Cal H}_3$ $|$ $h_{LR} > h_\infty\}$.

\vskip.2cm
%


%




\ 

\noindent{\bf Case 1)-b:}  Suppose that $u_L \in {\Cal H}_1$ and $u_R\in {\Cal H}_1$.
Here one has, by definition
$$
h_{LR} = {1\over 2c_1} (c_1 w_L + c_1 w_R + v_R-v_L) = {w_L + w_R
\over 2} + {v_R - v_L\over 2c_1}.
\leqno{(2.15)}
$$
According to our formulation (1.16) of the initiation criterion, the
solution will take its values in phase~1 region only in case  $h_{LR}
\le w_M$, while a phase~3 state will appear in the solution if $h_{LR}$
exceeds $w_M$. We distinguish between these two situations.

\

\noindent{\bf Case 1)-b1:} Suppose moreover that $-1 < h_{LR} \le w_M$.

We seek the solution in the form of three constant states separated
by a $-c_1$ contact wave and a $c_1$ contact wave:
$$
u(t,x) = \cases
 u_L \hfil & \text{ for } \quad x< -c_1t, \cr
u_* \hfil & \text{ for }\quad -c_1 t < x < c_1t, \cr
u_R \hfil & \text{ for } \quad x>c_1t.\cr
\endcases
\leqno{(2.16)}
$$
The intermediate state $u_* = (v_*,w_*)\in {\Cal H}_1$ must satisfy
the jump conditions:
$$
-c_1 (w_* - w_L) + v_* - v_L = 0, \qquad
c_1(w_R-w_*) + v_R - v_* = 0,
$$
which lead us to explicit expressions
$$
v_* = v_L + c_1(h_{LR} - w_L) \quad\hbox{ and }\quad w_* = h_{LR}.
\leqno{(2.17)}
$$
Because of the assumption $-1 < h_{LR} \le w_M$, it is immediate
that $w_*$ belongs to ${\Cal H}_1$. Note that no solution taken its values in $\Cal H$ exists when $h_{LR}<-1$.

For further reference, we state:

\

\noindent{\bf Lemma 2.3.} {\it Consider the Riemann problem} (1.1), (2.1) {\it in
case} 1)-b1 {\it i.e. when} $u_L\in {\Cal H}_1$, $u_R\in {\Cal H}_3$
{\it and} $-1 < h_{LR}\le w_M$. {\it Then the admissible weak solution
to these problems is given by} (2.16), (2.17). {\it The state}
$u_* = (v_*, w_*)$ {\it is a} ${\Cal C}^\infty$ {\it function of the
initial states} $u_L$ {\it and} $u_R$. {\it Moreover, when} $h_{LR}$
{\it tends to} $w_M$, {\it one has:}
$$
\lim_{h_{LR}\to w_M^-} v_* = v_L + c_1( w_M - w_L), \qquad
\lim_{h_{LR}\to w_M^-} w_* = w _M, 
\leqno{(2.18a)}
$$
and
$$
\lim_{h_{LR}\to w_M^-} {\partial v_*\over\partial h_{LR}} (u_L,
h_{LR}) = c_1, \qquad \lim_{h_{LR}\to w_M^-} {dw_*\over dh_{LR}}=1.
\leqno{(2.18b)}
$$

\vskip.2cm

\noindent{\bf Remark 2.3.} 1)~It is of interest to note that the assumption $-1 < h_{LR}\le w_M$
in Lemma~2.3 is always fulfilled if both $u_L$ and $u_R$ belong to ${\Cal H}_1$
and $|u_R-u_L|$ is small enough. This is clear in view of (2.15).

\noindent 2)~
In [2], an initiation criterion was introduced in
the case 1)--b1. Indeed, instead of the solution (2.16) containing no
phase boundary, the criterion in [2] selects in some cases a solution containing two phase boundaries (cf.~(2.19) below).

\

\noindent {\bf Case 1)-b2:}  Suppose now that $h_{LR} > w_M$.

According to our initiation criterion (1.16), the solution must contain
(at least) one phase~3 state.  We seek the solution in the form
$$
u(t,x) =  \cases
u_L \hfil & \text{ for } \quad x< -c_1t,\cr
u_1 \hfil & \text{ for } \quad -c_1 t < x < V', \cr
u_2 \hfil & \text{ for } \quad V'< x< V,\cr
u_3 \hfil & \text{ for } \quad V<x<c_1t,\cr
u_R \hfil & \text{ for } \quad x>c_1t, \cr
\endcases 
\leqno{(2.19)}
$$
where $u_1, u_3 \in {\Cal H}_1$ and $u_2\in {\Cal H}_3$ and $-c_3 < V'
< 0 < V < c_3$.  The jump conditions read:
$$
\left\{
\eqalign{
-c_1(w_1-w_L) + v_1 - v_L = 0, \qquad
&c_1(w_R-w_3) + v_R-v_3 = 0, \cr
V'(w_2-w_1)+ v_2 -v_1 = 0, \qquad
&V'(v_2-v_1)+c_3^2w_2-c_1^2w_1 = 0, \cr
V(w_R-w_3) + v_R - v_3 = 0, \qquad
&V(v_R-v_3) + c_3^2w_R-c_1^2w_3 = 0. \cr
}
\right.
\leqno{(2.20a)}
$$
They must be completed by the kinetic relation along the lines $x/t=V'$
and $x/t=V$:
$$
\left\{
\eqalign{
&{1\over 2} (k_1-k_3) (w_Mw_m - w_1 w_2) = -\phi (-V')\cr
&{1\over 2} (k_1 - k_3) (w_M w_m - w_2 w_3) = \phantom{-}\phi(V). \cr
}
\right.
\leqno{(2.20b)}
$$
It can be shown that, in fact, $V'=V$ (cf.~our previous calculations,
(2.7)). As for case~1)-a1, it can be checked that (2.20) determine
uniquely the admissible solution. We omit the details
and simply state our result of continuous dependence.

\

\noindent{\bf Lemma 2.4.} {\it Consider the Riemann problem} (1.1), (2.1) {\it in 
case} 1)-b2, {\it i.e. when} $u_L\in {\Cal H}_1$, $u_R\in {\Cal H}_1$ 
{\it and} $h_{LR} > w_M$. {\it Then the admissible weak solution of this 
problem is given by} (2.19), (2.20). {\it One can consider the states} 
$u_1$, $u_2$, $u_3$ {\it and the speed} $V$ {\it as functions of} $u_L$ 
{\it and} $u_R$.  {\it Then the functions} $u_i$ {\it and} $V$ {\it are} 
${\Cal C}^1$ {\it (with Lipschitz coninuous derivatives) and one has}
$$
\lim_{h_{LR}\to w_M} V = 0, \qquad \lim_{h_{LR}\to w_M} v_j = v_L
+ c_1(w_M-w_L) \quad\hbox{ for }\quad j=1,2 \hbox{ or } 3, 
\leqno{(2.21a)}
$$
$$
\lim_{h_{LR}\to w_M} w_1 = \lim_{h_{LR}\to w_M} w_3 = w_M,\qquad
\lim_{h_{LR}\to w_M} w_2 = {k_1\over k_3} \, w_M,
\leqno{(2.21b)}
$$
and
$$
\lim_{h_{LR}\to w_M } {\partial v_j \over \partial h_{LR}} = c_1, \qquad
\lim_{h_{LR}\to w_M } {\partial w_j\over\partial h_{LR}} = 1 \quad\hbox{ for }
\quad j=1 \hbox{ or } 3.
\leqno{(2.21c)}
$$

\vskip.2cm

Note that the limits found in (2.18) and (2.21) for the functions $v_*$
and $w_*$ and $v_j$ and $w_j$ (for $j=1$ or 3) coincide. Hence if in 
case 1)-b1, we set 
$$
v_j = v_*, \qquad w_j = w_* \quad\hbox{ for }\quad j=1 \hbox{ or } 3,
\leqno{(2.22)}
$$
then the functions $v_j$ and $w_j$ are globally of class ${\Cal C}^1$
with Lipschitz continuous derivatives in the whole range of values
$\{u_L\in {\Cal H}_1$, $u_R\in {\Cal H}_3\}$.

From Lemmas~2.1 to 2.4, we deduce the following property of
{\it continuous dependence} of the solution of the Riemann problem.

\

\noindent{\bf Theorem 2.1.} {\it Consider the admissible weak solution to the
Riemann problem} (1.1), (2.1) {\it described in Lemmas}~2.1 {\it to} 2.4.
{\it Then the states and the wave speeds in the solution are locally Lipschitz
continuous functions of the initial constant states} $u_L$ {\it and} $u_R$.
{\it As a consequence, if} $u_1(\cdot,0)$ {\it and} $u_2(\cdot,0)$
{\it are two Riemann initial data for system} (1.1), {\it the corresponding
admissible solutions} $u_1$ {\it and} $u_2$ {\it satisfy the following}
$L^1$ {\it continuous property:}
$$
\int_A^B |u_2 (t,x) - u_1(t,x)|\, dx \le O(1)\int_{A-c_1t}^{A+c_1t}
|u_2(0,x) - u_1(0,x)|\, dx
\leqno{(2.23)}
$$
{\it for all} $A < B$ {\it and} $t\ge 0$.

\vskip.2cm

\

We now turn to the Riemann problem in the half space $]0,\infty[$, i.e.
problem (1.1), (2.2). Two cases must be distinguished:

Case~2)--a: $\quad u_0\in {\Cal H}_3$,

Case~2)--b: $\quad u_0\in {\Cal H}_1$.

\vskip.2cm

\noindent {\bf Case 2)--a:}  Suppose that $u_0\in {\Cal H}_3$.

\vskip.2cm

According to condition (1.20)$_i$, the solution must contain a phase
boundary if and only if $h_0 < w_m^{cr}$ where
$$
h_0 = w_0 + {v_0\over c_3}
\leqno{(2.24)}
$$
and $w_m^{cr}$ is the critical value for initiation introduced in Section~1.
We recall that $w_0\ge w_m$ and $w_m^{cr}\ge w_m$.

\

\noindent{\bf Case~2)--a1:}  Suppose moreover that $h_0\ge w_m^{cr}$.

Then the solution $u$ to (1.1), (2.2) must stay entirely in phase~3,
so we seek $u$ in the form
$$
u(t,x) =\cases
u_* \hfil & \text{ for }\quad x<c_3t,\cr 
u_0\hfil & \text{ for }\quad x>c_3t.\cr
\endcases
\leqno{(2.25)}
$$
To satisfy the boundary condition (1.17a), we must have
$$
v_* = 0.
\leqno{(2.26a)}
$$
From the Rankine-Hugoniot relation
$
c_3 (w_0-w_*) + v_0 - v_* = 0,
$
we deduce $w_*$ :
$$
w_* = w_0 + {v_0\over c_3} = h_0.
\leqno{(2.26b)}
$$
The solution is completely determined by (2.25), (2.26).

\

\noindent{\bf Case~2)--a1:}  Suppose now that $h_0 < w_m^{cr}$.

In that case, the solution must contain a phase boundary, so we set
$$
u(t,x) = \cases 
u_-\hfil & \text{ for }\quad x<Vt, \cr 
u_+ \hfil& \text{ for }\quad Vt < x < c_3t, 
\cr u_0\hfil& \text{ for }\quad x>c_3t\cr
\endcases
\leqno{(2.27)}
$$
In view of the boundary condition (1.17a), one has
$$
v_- = 0.
$$
We determine $V$, $w_-$ and $u_+ = (v_+, w_+)$ by writing the 
Rankine-Hugoniot relations satisfied along the lines $x/t=V$ and
$x/t=c_1$, as well as the kinetic relation along $x/t=V$. By calculations 
similar to those made in case~1)--a1, we obtain the following formulas:
$$
v_- = 0, \qquad w_- = {c_3+V\over c_1+V} \, h_0,
\leqno{(2.28a)}
$$
$$
v_+ = v_0 + c_3 w_0 - {c_1-V\over c_3-V} \, c_3 h_0 = -{c_1-c_3
\over c_3-V}\, c_3h_0, \quad w_+ = {c_1-V\over c_3-V} \, h_0,
\leqno{(2.28b)}
$$
where $V$ is given by the following implicit equation:
$$
{k_1-k_3\over 2} \big\{ w_Mw_m - {(c_3+V)(c_1-V)\over (c_1+V)(c_3-V)}
\, h_0^2\big\} = \phi(V).
\leqno{(2.28c)}
$$
These formulas determine the solution in this case.
%




\

\noindent {\bf Case 2)--b:}  Suppose that $u_0\in {\Cal H}_1$:

According to condition (1.20)$_{ii}$, we have to distinguish between two 
cases. Here $h_0 = w_0 + \displaystyle{v_0\over c_1}$.

\

\noindent{\bf Case 2)--b1:}  Suppose moreover that $h_0\le w_M^{cr}$.

Then the solution stays entirely in phase~1:
$$
u(t,x) = \cases 
u_* \hfil & \text{ for }\quad x<c_1t,\cr
u_0\hfil & \text{ for }\quad x>c_1t \cr
\endcases
\leqno{(2.29)}
$$
where
$$
v_* =0 \quad \hbox{ and } \quad w_* ={v_0\over c_1} + w_0 = h_0\in{\Cal H}_1.
\leqno{(2.30)}
$$

\noindent{\bf Case 2)--b2:}   Suppose moreover that $h_0 > w_M^{cr}$.

Then the solution contains a phase boundary, i.e.
$$
u(t,x) = \cases 
u_- \hfil & \text{ for } \quad x<Vt\cr
u_+ \hfil & \text{ for } \quad Vt < x < c_1t \cr
u_0 \hfil & \text{ for } \quad x >c_1t, \cr
\endcases 
\leqno{(2.31)}
$$
where $u_- \in {\Cal H}_3$, $u_+ \in {\Cal H}_1$ and $V \in\, ]0,c_1[$.
The states $u_-$, $u_+$ and the speed $V$ are uniquely determined by (1.17a),
the Rankine-Hugoniot relations and the kinetic relation. We omit the details.

Finally, we conclude with the result of $L^1$ continuous dependence for
the Riemann problem in a half space.

\

\noindent{\bf Theorem 2.2.} {\it Consider the admissible weak solution of the 
Riemann problem} (1.1), (2.2) {\it described by cases}~2). {\it The
states and the wave speeds in the solution are Lipschitz continuous
functions of the initial state} $u_0$. {\it As a consequence, if} $u_0'$
{\it and} $u''_0$ {\it are two Riemann data for the system} (1.1) {\it in the
half space, then the corresponding admissible solutions} $u'$ {\it and}
$u''$ {\it satisfy}
$$
\int_A^B |u'(t,x) - u''(t,x)| \, dx \le O(1) (B+c_1t -\max (0,A-c_1t)\, )
|u'_0-u''_0|
$$
{\it for all} $0 < A< B$ {\it and} $t\ge 0$.

\vskip.2cm

\vfill

\

\noindent{\bf 3. EXISTENCE VIA GLIMM'S SCHEME: STABILITY}

\vskip.2cm

This section and the following one deal with the application of the random choice
 method, introduced by Glimm [15] for hyperbolic problems, to the system
of mixed type (1.1).  Our main result establishes existence of
a class of admissible weak solutions to (1.1). This serves to justify
the formulation of the Cauchy problem proposed in Section~1.

Based on successive resolutions of Riemann problems, Glimm's method yields
a sequence of approximate solutions for the Cauchy problem associated
with (1.1). Our goal is to prove the convergence of these approximate
solutions to an admissible weak solution to the problem, in case
the initial data is a small BV perturbation of a single propagating phase
boundary. The main result of this section, Theorem~3.1, yields the
stability of the scheme in the BV norm. This guarantees its convergence
in the $L^1$ norm to a function of bounded variation, which indeed
is a weak solution of (1.1). Showing that this function is an
admissible solution requires a more detailed analysis which is performed
in the next section.

It is emphasized that a phase boundary is a wave with (necessarily)
{\it large} strength. Our result of stability here is related
to the ones obtained by Chern [5] and Schochet [37] who treated Glimm's
scheme with  large data for strictly hyperbolic systems.

We consider the system (1.1) on the whole line ($x\in \RR$) with the following
Cauchy data:
$$
u(0,x) = u_0(x) = \cases 
u_L^0 (x) = (v_L^0 (x), w_L^0(x))\hfil
& \text{ for }\quad x<0,\cr
\cr
u_R^0 (x) = (v_R^0(x),w_R^0(x))\hfil & \text{ for }\quad x>0.\cr
\endcases 
\leqno{(3.1a)}
$$
The functions $u_L^0 \in BV_{loc}(\RR_-, {\Cal H})$ and $u_R^0 \in BV_{loc}
(\RR_+,{\Cal H})$ are assumed to be close to two given constant
states $u_L^* = (v_L^*, w_L^*)$ and $u_R^* = (v_R^*, w_R^*)$ respectively, i.e.
$$
\Vert u_L^0 - u^*_L\Vert_{BV(\RR_-)} + \Vert u_R^0 - u^*_R\Vert_{BV(\RR_+)}
\ll 1.
\leqno{(3.1b)}
$$
For definiteness, we consider the case that $u_L^* \in {\Cal H}_1$ 
and $u_R^*\in {\Cal H}_3$. We are assuming that the Riemann problem with data
$u_L^*$ and $u_R^*$ is solved by a unique wave, i.e. a single phase boundary
but not contact discontinuity. This assumption allows us to focus our
attention on phase boundaries which are the main difficulty regarding 
system (1.1). Let $u^*$ be the solution of this Riemann problem; for some
speed $V^*$, one has
$$
u^*(t,x) = \cases
 u_L^* \hfil& \text{ for }\quad x<V^*t,\cr
\cr
u_R^*\hfil& \text{ for }\quad
x>V^*t.\cr
\endcases
\leqno{(3.2a)}
$$
In the case of a characteristic phase boundary, i.e. when $V^*=c_3$,
we will restrict ourselves to the case that {\it no strong wave\/} arises from a 
perturbation of the initial states $u_L^*$ and $u_R^*$. According to 
Lemma~2.1 of Section~2 (cf.~formulas (2.10a) and (2.10b)), this is the case under the following assumptiom:
condition:
$$
\hbox{If } V^* =c_3, \hbox{ then } 
w_L^* = 0 \hbox{ and } w_R^* = \sqrt{ {\phi'(c_3)\over c_3} }.
\leqno{(3.2b)}
$$
Note that (3.2b) implies that $h_{LR}=0$ and $v_R=v_L-\sqrt{c_3\phi'(c_3)}$. So a Riemann problem with initial data in the neighborhood of $u_L$ and $u_R$ takes its values in the same neighborhood. 

We shall prove that problem (1.1), (3.1) admits an admissible weak
solution which has the following structure:
$$
u(t,x) =\cases
 u_L(t,x) \hfil & \text{ for } \quad x < \chi(t), \cr
\cr
u_R(t,x)\hfil & \text{ for } \quad x>\chi(t), \cr
\endcases
\leqno{(3.3a)}
$$
where
$$
u_L\in L^{\infty}_{loc}([0,\infty[, BV(\RR,{\Cal H}_1)) \hbox{ and }u_R\in 
L^{\infty}_{loc}([0,\infty[, BV(\RR,{\Cal H}_3))
\leqno{(3.3b)}
$$
and
$$
\chi \in W_{loc}^{1,\infty} ([0,\infty[,\RR) \quad\hbox{ and }\quad
{d\chi \over dt} \in BV_{loc}([0,\infty[,\RR).
\leqno{(3.3c)}
$$
Setting
$$
\tilde u^*(t,x) = \cases
 u_L^* \hfil & \text{ for } \quad x < \chi(t),\cr
u_R^* \hfil & \text{ for } \quad x>\chi(t),\cr
\endcases 
$$
we shall also show that for all times $T>0$
$$
\Vert u(T) -\tilde u^*(T)\Vert_{L^\infty(\RR,{\Cal H})} + TV_{\RR}
(u(T)-\tilde u^*(T) ) + TV_0^T ({d\chi_0\over dt}-V^*)\ll 1.
\leqno{(3.4)}
$$
The solution will be obtained as the limit of approximate solutions to
(1.1), (3.1), having a structure similar to the one described by (3.3), (3.4).
These approximate solutions are given by Glimm's scheme which we now describe.

\

\

Let $\tau > 0$ and $h>0$ be time and space mesh sizes satisfying
the CFL condition $\tau c_1 < h$. The ratio $\lambda = h/\tau$ is taken
to be a constant. Let $\{ a_n\}_{n\ge 1}$ be an equidistributed sequence
with values in the interval $]-1,1[$. We define $u^h(0,x)$ by 
$L^2$--projection from the data $u_0$:
$$
u^h(0,x) = {1\over 2h} \int_{mh}^{(m+2)h} u_0(y) \, dy \quad\hbox{ for }\quad
x\in [m h, (m+2)h[ \hbox{ with } m \hbox{ even.}
\leqno{(3.5a)}
$$
Note here that $|u^h (0,x)-u_L^*|\ll 1$ for $x<0$ and $|u^h(0,x)-u_R^*|
\ll 1$ for $x>0$; also $u^h$ satisfies $u^h(0,x)\in {\Cal H}_1\cup {\Cal H}_3$.
If $u^h$ is known up to the time $t=n\tau -0$, we define $u^h(n \tau +0, x)$
by a random choice projection using $a_n$:
$$
u^h(n\tau +0,x) = u^h(n\tau -0, (m+1+a_n)h-0) 
\leqno{(3.5b)}
$$
for $x\in [mh, (m+2)h[$ with $m+n$ even. Then the approximate solution $u^h$
in the strip $\{ n\tau \le t < (n+1)\tau\}$ is computed by solving
the Riemann problems for system (1.1) at each center $x=mh$ with $m+n$ even.

As a consequence of our result of stability below, this construction indeed
makes sense and yields $u^h(t,x)$ for all times $t\ge 0$. In particular,
because of the assumption (3.2b),
the values $u^h(t,x)$ stay in the neighborhoods of $u_L^*$ or $u^*_R$.
This implies that case 1)-b2 of Section~2 will never occur here. The
possible wave structures of the Riemann problem used in the construction of
$u^h$ can be listed.

\vskip.2cm

\noindent{\bf Remark 3.1.} As a very first step toward a general proof of convergence
of $u^h$, we may consider the case when $u_L^0$ and $u_R^0$ are constant,
equal to $u_L^*$ and $u_R^*$ respectively. In that case, $u^h$ can
be computed explicitly and consists for each time $t$ of a single phase
discontinuity connecting $u_L^*$ at the left to $u_R^*$ at the right.
The position of the phase discontinuity, say $\chi^h(t)$, is shifted
to the left or to the right (depending on $a_n$ and the speed $V^*$)
at each time $n\tau$.  This is typical behavior for Glimm's scheme, which,
as is well known, does not produce any numerical diffusion of the
discontinuities. Using only the equidistributedness of $\{a_n\}$, it is
an easy matter to show that 
$\chi^h(t)$ converges to $V^* t \equiv \chi(t)$ for each time $t\ge 0$.

\vskip.2cm

According to the technique of Glimm, the first step toward a proof of uniform
BV stability for the scheme consists of studying wave interactions
between Riemann wave patterns.  Let $R(u_\ell,u_r)$ be the solution of the Riemann problem with initial data $u_\ell$
at the left and $u_r$ at the right (cf. Section~2 for the explicit construction). The wave strengths are defined first
in case $V^*\ne c_3$.  We denote by ${\Cal E}_1(u_\ell,u_r)$
and ${\Cal E}_2(u_\ell,u_r)$ the strengths of the left and right contact
waves in $R(u_\ell, u_r)$ respectively. By convention ${\Cal E}_2(u_\ell,u_r)=0$
when $R(u_\ell,u_r)$ contains a supersonic phase boundary so that
no right contact discontinuity is present. We denote by ${\Cal E}_0(u_\ell,u_r)$
the strength of the phase boundary in $R(u_\ell,u_r)$  in case where
$u_\ell$ and $u_r$ are in diffferent phases.  By convention, the strengths are
always measured in terms of the jump of the variable $w$ across
the wave under consideration. When a phase boundary is present, we denote its speed 
by ${\Cal V}(u_\ell,u_r)$.

Consider now the case that $V^* = c_3$. We define ${\Cal E}_0$, ${\Cal E}_1$,
${\Cal E}_2$ and ${\Cal V}$ in the same way as above, except in case that 
the Riemann problem $R(u_\ell,u_r)$ admits a supersonic phase boundary.
In this latter case, we virtually split the phase discontinuity into two
distinct waves and set:
$$
{\Cal E}_2 (u_\ell, u_r) = w_r -w_+\big|_{h_{\ell r}=0} - h_{\ell r}
{\partial w_+\over\partial h_{\ell r}} \big|_{h_{\ell r}=0}  =
 w_r -\sqrt{ {\phi'(c_3)\over c_3} } - 2 c_3 h_{\ell r} / (c_1 + c_3),
$$
and
$$
{\Cal E}_0 (u_\ell, u_r) = w_+\big|_{h_{\ell r}=0} + h_{\ell r} 
{\partial w_+ \over \partial h_{\ell r}}\big|_{h_{\ell r}=0} - w_-  =  
\sqrt{ {\phi'(c_3)\over c_3} } + 2 c_3 h_{\ell r} / (c_1 + c_3) - w_-,
$$
where $w_-$ and $w_+$ are the values taken by the solution of $R(u_\ell,u_r)$
at the left and at the right of the phase boundary respectively and $h_{\ell r}$ is 
given by (2.3). We recall that
$\displaystyle{\partial w_+\over\partial h_{\ell r}}\big|_{h_{\ell r}=0}$ is given
 by Lemma~2.1. In other words, we extend the definition of ${\Cal E}_0$ and
 ${\Cal E}_2$, known for $h_{\ell r} >0$, to negative values of $h_{\ell r}$ 
so that their extensions are of class ${\Cal C}^1$.

From the results in Lemmas~2.1 and 2.2, one easily checks that
$$
\left\{
\eqalign{
&\hbox{the functions } {\Cal E}_1, {\Cal E}_2, {\Cal E}_0 \hbox{ and }
{\Cal V} \hbox{ are Lipschitz continuous}\cr
&\hbox{functions of their arguments } (u_\ell,u_r),
\hbox{ the }  {\Cal E}_i \hbox{ 's are of class}\cr
&{\Cal C}^1 \hbox{ (with Lipschitz continuous derivatives) away from } {\Cal V} = 0. \cr
}
\right.
\leqno{(3.6)}
$$
In (3.6) and all the sequel, only states which are close to either
$u_L^*$ or $u_R^*$ are considered.

The wave interaction estimates are derived in the following lemma:

\vskip.2cm

\noindent{\bf Lemma 3.1.} {\it Consider states} $u_\ell$, $u_p$ {\it and}
$u_r$ {\it which are close to either} $u_L^*$ {\it or} $u_R^*$.

1) {\it If} $u_\ell \in {\Cal H}_1$, $u_p \in {\Cal H}_1$
{\it and} $u_r \in {\Cal H}_3$, {\it then for} $j=0,1,2$ {\it we have}
$$
{\Cal E}_j (u_\ell,u_r) = {\Cal E}_j (u_\ell,u_p) + {\Cal E}_j(u_p,u_r)
+ O(1)|{\Cal E}_2(u_\ell,u_p) |
\leqno{(3.7a)}
$$
{\it and}
$$
{\Cal V}(u_\ell,u_r) = {\Cal V} (u_p,u_r) + O(1)|{\Cal E}_2 (u_\ell,u_p)|.
\leqno{(3.7b)}
$$

2) {\it If} $u_\ell \in {\Cal H}_1$, $u_p\in {\Cal H}_3$ {\it and}
$u_r\in {\Cal H}_3$, {\it then for} $j=0,1,2$ {\it we have}
$$
{\Cal E}_j (u_\ell,u_r) = {\Cal E}_j (u_\ell, u_p) + {\Cal E}_j (u_p,u_r)
+ O(1) \, D(u_\ell, u_p,u_r)
\leqno{(3.8a)}
$$
{\it and}
$$
{\Cal V} (u_\ell,u_r) = {\Cal V} (u_\ell, u_p) + O(1)\, \bigg\{|{\Cal E}_1 (u_p, u_r) | + |{\Cal E}_2 (u_p, u_r) |\bigg\},
\leqno{(3.8b)}
$$
{\it where we have set}
$$
D(u_\ell,u_p,u_r) = \cases
|{\Cal E}_1 (u_p, u_r) |\hfil & \text{ if }\quad {\Cal V}(u_\ell,u_p)\le c_3, \cr
\cr
|{\Cal E}_1 (u_p,u_r) | + ({\Cal V}(u_\ell,u_p)-c_3)|{\Cal E}_2(u_p,u_r)|
\hfil & \text{ if } \quad {\Cal V}(u_\ell,u_p) \ge c_3. \cr
\endcases 
\leqno{(3.9)}
$$

3) {\it if either} $u_\ell$, $u_p$ {\it and} $u_r \in {\Cal H}_1$
{\it or} $u_\ell, u_p$ {\it and} $u_r\in {\Cal H}_3$, {\it then for}
$j=1,2$ {\it we have}
$$
{\Cal E}_j (u_\ell,u_r) = {\Cal E}_j (u_\ell,u_p) + {\Cal E}_j (u_p,u_r).
\leqno{(3.10)}
$$

\vskip.2cm

\noindent{\bf Remark 3.2.} 1) Estimates in Lemma~3.1 mainly contain {\it linear} interaction terms 
instead of quadratic ones as is the case in [15]. Linear error
terms were previously found useful by [5] and [37] to treat {\it strictly
hyperbolic} systems with {\it large} data. In (3.9), the interaction term
is proportional to the {\it angle} between the $c_3$ contact discontinuity
and the phase boundary. Such a term was used by Liu [34] to analyze
{\it non-genuinely nonlinear} systems of conservation laws.

\noindent 2) When $V^*\ne c_3$, the derivation of the estimates in Lemma~3.1 requires 
only the Lipschitz continuity of the functions ${\Cal E}_j$ and ${\Cal V}$, 
which is exactly the regularity available in general.

\noindent 3) The smallness condition (on $|u_{\ell} - u^*_L|$, etc) in Lemma~3.1
is necessary only to prevent initiation of a new phase when solving a Riemann
problem with data in a single phase.

\vskip.2cm

\

%






 
\noindent{\bf Proof of Lemma 3.1.} We first give the proof of (3.7), and then that of (3.8).
The proof of (3.10) is trivial.

In view of the results of Section~2, the functions ${\Cal E}_0$, ${\Cal E}_1$, ${\Cal E}_2$
and ${\Cal V}$ are (at least) Lipschitz continuous functions of their arguments.
Hence the formulas (3.7) will follow easily if we check that
$$
{\Cal E}_j (u_\ell, u_r) = {\Cal E}_j (u_\ell,u_p) + {\Cal E}_j(u_p,u_r)
\leqno{(3.11a)}
$$
and
$$
{\Cal V} (u_\ell,u_r) = {\Cal V} (u_p,u_r)
\leqno{(3.11b)}
$$
hold whenever ${\Cal E}_2 (u_\ell,u_p)=0$, i.e. when there is no
right wave in the left wave packet $R(u_\ell,u_p)$.  But this last
statement is obvious because the left-waves in $R(u_\ell,u_p)$ and
$R(u_p,u_r)$ are associated with a linearly degenerate characteristic
field.  Such waves can be superimposed without any interaction and the 
wave strengths are simply summed up, cf.~(3.11a).  The speed of the phase
boundary remains unchanged, cf.~(3.11b). (These facts can be
checked directly from the analytical expressions in Section~2.) The proof
of (3.7) is completed.

We now prove (3.8). We notice first that (3.11a) as well as 
$$
{\Cal V}(u_\ell,u_r) ={\Cal V}(u_\ell,u_p)
\leqno{(3.11c)}
$$
do hold provided that $D$ given by (3.9) vanishes. 

Namely, when ${\Cal V}(u_\ell,u_p)\le c_3$ and if $D(u_\ell,u_p,u_r)=0$, that
means that the right wave packet does not contain a $-c_3$ contact wave.
In that situation, the two wave patterns can be superimposed, without
any interaction. When ${\Cal V}(u_\ell,u_p) \ge c_3$ and if 
$D(u_\ell,u_p,u_r)=0$,
then the right wave packet does not contain a $-c_3$ contact wave and,
moreover, either it also has no $c_3$ contact wave or the speed 
${\Cal V}(u_\ell,u_p)$ equals $c_3$. In both situations, the left and right
wave patterns can be superimposed. Again, there is not interaction.
This proves (3.11a) and (3.11c).

When ${\Cal V}(u_\ell,
u_p)\le c_3$, estimates (3.8) follow from (3.11a), (3.11c) and the Lipschitz
continuity of ${\Cal E}_j$ and ${\Cal V}$. When ${\Cal V}(u_\ell,u_p)>c_3$,
since ${\Cal V}$ is away
from ${\Cal V}=0$ and according to the results in Section~2, the functions
${\Cal E}_j$ and ${\Cal V}$
are of class $W^{2,\infty}$. This allows us to apply the classical lemma
of division (e.g. [Ho]) and again to deduce (3.8) from (3.11a), (3.11c). 
The proof of (3.8) is completed. ~~$\Box$

\vskip.2cm

We will use the technique of Glimm to deduce from Lemma~3.1 the result
of BV stability of the scheme. We refer to [15] and [16] for the terminology
we use here. At this stage, we have to define functionals to control the total variation of
the solutions. The choice we propose here is motivated by the form of the terms of 
interaction found in Lemma~3.1. Note that the phase boundary which is a
``strong wave'' are treated separately from the ``small waves''.

The $(t,x)$-plane is divided into a set of diamonds $\Delta_{m,n}$
with centers $(n\tau,mh)$ ($n+m$ even) and with vertices
$$
\eqalign{
N= (\, (n+1)\tau, (m+a_{n+1})h), \quad 
&E= (n\tau , (m+1 + a_n)h), \cr
W= (n\tau, (m-1+a_n)h), \qquad \quad
&S= (\, (n-1)\tau, \, (m+a_{n-1})h).\cr
}
$$
Given a diamond $\Delta_{mn}$, we denote by $u_N$, $u_E$, $u_W$ and $u_S$
the values taken by $u^h$ at the vertices $N$, $E$, $W$ and $S$ respectively.

We give now the definition of the {\it approximate phase boundary}
in $u^h$ that we denote by $\chi^h\colon \, \RR_+\to \RR$. First of all,
it is a simple (but useful) observation that the phase boundary in $u^h$
is actually located at a single space position for each time $t=n\tau$. In other
 words, there is no spreading of the phase boundary.
Let $\chi^h$ be the piecewise linear curve which is discontinuous at each 
$t=n\tau$ and coincides with the phase boundary in $u^h$ inside each slab
$[n\tau , (n+1)\tau[$. Let ${\Cal D}$ be the set of all diamonds that
are crossed out by the phase boundary $\chi^h$.

We then introduce several functionals defined on space-like curves, say $J$,
passing through vertices of diamonds. Define
$$
L(J) = \sum \big( |{\Cal E}_1| + |{\Cal E}_2|\big)
\leqno{(3.12a)}
$$
the summation being on all small waves crossing the curve $J$, and
$$
B(J) = |{\Cal E}_0|
\leqno{(3.12b)}
$$
where ${\Cal E}_0$ is the strength of the phase boundary when crossing
the curve $J$.  The functional $L(J)$ bounds the total
variation of $u^h$ along the curve $J$ on both sides of the phase
boundary. $B(J)$ measures the jump of $u^h$ across the phase
boundary.
Next we define the potential interaction $Q(\Delta)$ in a diamond $\Delta$
by
$$
Q(\Delta) =\cases
 0 \hfil & \text{ if } \quad \Delta \not\in {\Cal D},\cr
|{\Cal E}_2(u_W,u_S)| \hfil & \text{ if }\quad \Delta = \Delta_{m,n}\in {\Cal D}
\hbox{ and } mh < \chi^h (n\tau),\cr
|{\Cal E}_1 (u_S,u_E)| \cr
+\theta ({\Cal V}(u_W,u_S)-c_3)|{\Cal E}_2(u_S,u_E)| \cr
~ \hfil & \text{ if } \quad \Delta = \Delta_{m,n} \in {\Cal D} \hbox{ and }
mh\ge \chi^h(n\tau), \cr
\endcases
\leqno{(3.13)}
$$
where $\theta\colon \, \RR\to\RR$ is the function defined by
$$
\theta (y) = 0 \hbox{ for } y < 0, \,\, y \hbox{ for } y \ge  0.
$$
Finally the potential wave interaction $Q(J)$ associated with a curve $J$ is
$$ 
Q(J) = \sum_{waves ~ at ~ the \atop left ~ of ~ \chi^h} |{\Cal E}_2| +
\sum_{waves ~ at ~ the \atop right ~ of ~ \chi^h} \bigg\{ |{\Cal E}_1|
+ \theta (V-c_3) |{\Cal  E}_2|\bigg\},
\leqno{(3.14)}
$$
where $V$ is the speed of the phase boundary when crossing $J$ and the
summation being on all waves crossing $J$. Note that $Q(J)$ is a linear
functional in terms of wave strengths.

\vskip.2cm

\noindent{\bf Lemma 3.2.} {\it Let $K$ be a sufficiently large constant. Let $J_1$
and $J_2$ be two space-like curves, $J_2$ being a successor of $J_1$.
Then one has}
$$
L(J_2) + KQ(J_2) \le L(J_1) + KQ(J_1)
\leqno{(3.15a)}
$$
{\it and}
$$
B(J_2) + KQ(J_2)\le B(J_2) + KQ(J_2).
\leqno{(3.15b)}
$$

\vskip.2cm

\noindent{\bf Proof.} We need only prove (3.15) when $J_2$ is an immediate successor 
of $J_1$; the general case follows by induction. We check first the formula
$$
Q(J_2) - Q(J_1) \le -{1\over 2} Q(\Delta),
\leqno{(3.16)}
$$
where $\Delta$ is the diamond limited by $J_1$ and $J_2$.  If 
$\Delta \in {\Cal D}$ and the right wave packet contains the phase
discontinuity, then in view of (3.14), (3.13):
$$
Q(J_2) - Q(J_1) = -|{\Cal E}_2 (u_W,u_W)| = - Q(\Delta) \le -{1\over 2}
Q(\Delta).
$$
If $\Delta \in {\Cal D}$ and the left wave packet contains the phase 
discontinuity, then in view of (3.14)
$$
\eqalign{
Q(J_2)-Q(J_1) &= -|{\Cal E}_1 (u_S, u_E) | -\theta ({\Cal V} (u_W,u_S)-c_3)|
{\Cal E}_2 (u_S,u_E)| \cr
&\quad + \bigg\{\theta ({\Cal V} (u_W,u_E) - c_3)-\theta ({\Cal V}
(u_W,u_S) - c_3)\bigg\} \sum_{waves ~ on ~ the ~ right\atop side~of~\Delta}
|{\Cal E}_2| \cr
}
$$
Since $\theta$ is Lipschitz continuous and using definition (3.13) and
(3.8b), it follows that
$$
\eqalign{
Q(J_2) - Q(J_1) &= - Q(\Delta) + O(1) \, Q(\Delta) \sum_{waves ~ on ~the ~right
\atop side ~of~ \Delta} |{\Cal E}_2| \cr
&=Q(\Delta )\bigg\{ -1 + O(1) \, L(J_1)\bigg\} \le -{1\over 2} Q(\Delta)\cr
}
$$
where in the last inequality we have assumed that $O(1)\, L(J_1)\le {1\over 2}$.
The condition $L(J_1) \ll 1$ indeed is ensured by induction if
the initial total variation is small enough.
Let us content ourselves with checking that $L(J_0) \ll 1 $ where
$J_0$ is the curve connecting centers of diamonds on the lines $t=0$
and $t=t_1=\tau$. Namely, using the definition of the wave strengths and (3.2b), one gets
$$
L(J_0) = O(1)\bigg\{ TV_{-\infty}^0 (u_L^0) + TV_0^{+\infty}(u_R^0)
+\Vert u_L^0 - u^*_L\Vert _{L^\infty(\RR_-)} 
+\Vert u_R^0 - u^*_R\Vert _{L^\infty(\RR_+)} \bigg\}.
$$
The right hand side of the above formula is small in view of (3.1b).

If now $\Delta\not\in {\Cal D}$, then from (3.10) one has trivially 
$$
Q(J_2)-Q(J_1)\le 0 \quad\hbox{ and }\quad Q(\Delta) = 0.
$$
Henceforth, the proof of (3.16) is completed.

We next consider $L(J_2) -L(J_1)$. If $\Delta \in {\Cal D}$, then by (3.7a),
(3.8a), (3.12a), (3.13):
$$
\eqalign{
L(J_2) - L(J_1) &= |{\Cal E}_1 (u_W,u_E)| + |{\Cal E}_2(u_W,u_E)|\cr
&\quad - |{\Cal E}_1 (u_W,u_S)| - |{\Cal E}_2 (u_W,u_S)| - |{\Cal E}_1
(u_S,u_E)| - |{\Cal E}_2 (u_S,u_E)| \cr
&= O(1)\, Q(\Delta). \cr
}
$$
If $\Delta\not\in {\Cal D}$, one has
$$
L(J_2) = L(J_1) \quad\hbox{ and }\quad Q(\Delta) = 0.
$$
This proves the formula
$$
L(J_2) = L(J_1) + O(1)\, Q(\Delta).
\leqno{(3.17)}
$$
From (3.16) and (3.17), we easily deduce (3.15a) provided that the constant
$K$ in (3.15a) is large enough.

Finally, it can be proved similarly that 
$$
B(J_2) = B(J_1) + O(1)\, Q(\Delta)
\leqno{(3.18)}
$$
which implies (3.15b) in view of (3.16). ~~$\Box$

Lemma 3.2 provides a uniform bound for the total variation of $u^h(t)$ at times $t=t_n$. Since 
$$
TV(u^h(t)) \le O(1) L(J), \hbox{ for all times }t \in [t_n,t_{n+1}[,
$$
where $J$ is the curve lying between the lines $t=t_n$ and $t=t_{n+1}$, 
we obtain a uniform control of the total variation of $u^h$ for all times. 
Let us define the function $\tilde u^h\colon\,
\RR_+\times \RR\to {\Cal H}$ by:
$$
\tilde u^h(t,x) =\cases
u_L^* \hfil & \text{ if } \quad x < \chi^h(t), \cr
u_R^* \hfil & \text{ if } \quad x>\chi^h(t), \cr
\endcases
\leqno{(3.19)}
$$
where $\chi^h$ is the approximate phase boundary associated with $u^h$.
From Lemma~3.2, one deduces the following result of stability. 

\vskip.2cm

\noindent{\bf Theorem 3.1.} {\it The functions} $u^h$ {\it given by Glimm's scheme
applied to the mixed system} (1.1) {\it and the data} (3.1), (3.2) {\it satisfy the following stability estimates:}
$$
T V_{-\infty}^{+\infty} (u^h - \tilde u^h)(t) \le O(1) N_1, 
\leqno{(3.20a)}
$$
$$
\Vert u^h (t)-\tilde u^h(t)\Vert_{L^\infty (\RR, {\Cal H})} \le O(1) \big(N_1
+ N_2\big)
\leqno{(3.20b)}
$$
{\it and}
$$
\Vert u^h (t) - u^h(t')\Vert_{L^1(\RR, {\Cal H})} \le O(1)N_1\bigg( |t-t'|+h\bigg)
\leqno{(3.20c)}
$$
{\it for all times} $t\ge 0$ {\it and} $t'\ge 0$, {\it with} $N_1$ {\it and}
$N_2$ {\it defined by}
$$
N_1 = \cases
TV_{-\infty}^{0} (u_L^0) + TV_0^\infty (u_R^0)  \hfil & \text{ if } \quad V^*\ne c_3,\cr
\cr
TV_{-\infty}^0(u_L^0)+TV_0^\infty(u_R^0) + N_1 \hfil & \text{ if } \quad V^*=c_3 \cr
\endcases 
\leqno{(3.21a)}
$$
{\it and}
$$
N_2 = \Vert u_L^0 - u_L^* \Vert_{L^\infty (\RR_-, {\Cal H}_1)} +\Vert u_R^0
- u_R^*\Vert _{L^\infty(\RR_+,{\Cal H}_3)}.
\leqno{(3.21b)}
$$

\vskip.2cm

By Helly's theorem, the estimates (3.20) imply that (a subsequence of) 
$\{u^h\}$ converges in $L_{loc}^1$ strongly to a function $u$ as 
$h\to 0$.  This function has bounded variation in space and satisfies 
the same bounds as $u^h$ in (3.20). It is a classical matter 
(Glimm [15], Liu [33], [34]) to check that $u$ indeed is a {\it weak solution} 
to the system of conservation laws (1.1). It also satisfies the entropy 
inequality as well as the initial condition. It remains to show that 
$u$ is admissible, i.e.  satisfies the kinetic relation (cf.~Section~4).

\vskip.2cm

\noindent{\bf Remark 3.3.} If condition (3.2b) is violated, then perturbating
of a characteristic phase boundary produces a $c_3$--contact wave
with strong strength. Then, one would
have to deal with interaction between two strong waves traveling
with arbitrary close speeds. Initiation of new phases is possible. It is not clear whether the total
variation of $u^h$ would remain uniformly bounded in that case.

\vskip.2cm

\vfill
\eject
\noindent{\bf 4. EXISTENCE VIA GLIMM'S SCHEME: ADMISSIBILITY}

\vskip.2cm

In this section, we prove that the weak solution $u=\lim u^h$ found in Section~3
by Glimm's scheme does satisfy the kinetic relation (1.11). This establishes
that $u$ is an admissible weak solution to our problem and leads to the
desired result of existence and stability.

First of all, we notice that the kinetic relation (1.11) is formulated in
a pointwise sense, more precisely (1.11) must hold almost everywhere with 
respect to the {\it Hausdorff measure} $H_1$. However from the results
in Section~3, we only have that $u^h$ converges to $u$ at almost every point
with respect to the {\it Lebesgue measure} on $\RR_+\times \RR$. This
latter property is thus not sufficient to pass to the limit in the kinetic
relation.

We prove in this section a result of pointwise convergence for the 
approximate phase
boundary $\chi^h$, cf.~Theorem~4.1. This result is derived by using
the techniques introduced by Glimm and Lax in [16]. The focus of [16] 
was the case of a strictly hyperbolic system of two conservation laws 
with small data. Extensions of the results in [16] can also be found 
in the papers of Di~Perna [10] and Liu [34]. 
In our situation, we have a (special case of) a system
of mixed type with large data.

Next in Theorem 4.2, we prove that the above result is sufficient for the
passage to the limit in the kinetic relation, assuming that the speed
of the phase boundary $V^*$ does not vanish.

Let us consider the phase boundary $\chi^h\colon \, \RR_+\to\RR$ in the
approximate solution $u^h$. The function $\chi^h$ is discontinuous and
piecewise linear. It jumps up to a distance of $\pm 2h$ at each time step.
It is easy to verify the following lemma.

\vskip.2cm

\noindent{\bf Lemma 4.1.} {\it The function} $\chi^h\colon\, \RR_+\to\RR$ {\it satisfies
the following uniform estimate:}
$$
|\chi^h (t) -\chi^h(t')|\le {1\over\lambda} |t-t'| + 2h \quad\hbox{\it for }
\quad 0\le t\le t'.
\leqno{(4.1)}
$$

\vskip.2cm

By Ascoli's theorem, the sequence $\{\chi^h\}$ must converge
on each compact set in the uniform topology to a function $\chi\in
W_{loc}^{1,\infty}([0,\infty[,\RR)$.
The next lemma gives a bound for the total variation of the functions
$\dot \chi^h\colon \, \RR_+\to\RR$ defined by
$$
\dot \chi^h (t) = {d\chi^h(t) \over dt}\quad (\hbox{constant}) \hbox{ on each
interval } [n\tau, \, (n+1)\tau [.
\leqno{(4.2)}
$$
From an analysis of the waves crossing the phase boundary, we prove the
following result. (The proof is given after the statement of Theorem~4.1.)

\

\noindent{\bf Lemma 4.2.} {\it For all times} $T>0$, {\it one has the uniform estimate}
$$
TV_0^T(\dot\chi^h )\le O(1)\big\{ TV_{-(T/\lambda)-2h}^0 (v_L^0
- c_1 w_L^0) + TV_0^{(T/\lambda)+2h}(u_R^0)+N \big\},
\leqno{(4.3)}
$$
where $N=0$ if $V^*\ne c_3$ and $N=\Vert u_R^0-u^*_L\Vert _{L^\infty(
-(T/\lambda) - 2h)}+\Vert u_R^0-u^*_R\Vert _{L^\infty(0,
T/\lambda +2h)}$ if $V^*=c_3$.

\vskip.2cm

Hence, from Lemmas~4.1 and 4.2, the equidistributedness of the sequence 
$\{a_n\}$ and arguments in [16], we deduce the following pointwise 
convergence property.

\vskip.2cm

\noindent{\bf Theorem 4.1.} {\it The functions} $\chi^h$ {\it and} $\dot \chi^h$
{\it converge to the functions} $\chi$ {\it and} $\displaystyle{d\chi
\over dt}$ {\it respectively in the following sense:}
$$
\Vert\chi -\chi^h\Vert_{L^\infty(]0,T[,\RR)} \to 0\quad \hbox{\it when }\quad
h\to 0, \quad\hbox{\it for all }\quad  T>0
\leqno{(4.4)}
$$
{\it and}
$$
\dot \chi^h(t) \to {d\chi (t)\over dt} \quad\hbox{\it for all times }\quad
t \in \RR_+ \setminus E,
\leqno{(4.5)}
$$
{\it where} $E\subset \RR_+$ {\it is an at most countable set.}

\vskip.2cm

We give first the proof of Lemma~4.2 and then the one of Theorem~4.1.

\vskip.2cm

\noindent{\bf Proof of Lemma 4.2.} Let $T$ be fixed and let $N$ be such that 
$N\tau\le T < (N+1)\tau$. Recall that ${\Cal D}$ is the set of all
diamonds which contain a part of the phase boundary. Let $J$ be the
space-like curve which limits the domain of dependence of the
diamonds in ${\Cal D}$ with centers below the line $t=N\tau$.
For each time $t= n\tau$, $n=0,1,\ldots,N$, $J$ encloses a finite number
of diamonds that we denote by $\Delta _n^m$ for $m=1,2,\ldots,N+1-n$.
They are ordered increasingly. We define $m(n)$ to be such that
$\Delta_n^{m(n)} \in {\Cal D}$. 

By Lemma~3.1, the speed of $\chi^h$ at the time $n\tau$ is estimated from
the its value at time $(n-1)\tau$:
$$
\dot \chi (n\tau + 0) =\dot \chi ^h(\, (n-1)\tau +0) +O(1) 
|{\Cal E}(\Delta_{n,m(n)})|,
\leqno{(4.6)}
$$
where $|{\Cal E}(\Delta_{n,m(n)})|$ represent the strength of the waves entering the diamond.
By summation with respect to $n=1,\ldots,N$, we obtain
$$
\eqalign{
TV_0^{(N+1)\tau} (\dot \chi^h) &= \sum_{n=1}^N |\dot \chi^h(n\tau + 0)-
\dot \chi(\, (n-1) \tau +0)|\cr
&=O(1)\sum_{n=0}^{N-1} |{\Cal E}(\Delta_{n,m(n)})|, \cr
}
$$
which bounded by the total variation of $u^h$ measured along both sides of
 the phase boundary. By using conservation laws for wave stregnths, as made 
in [16], one could check that the total variation of $u^h$ along this curve 
is bounded by the initial total variation. Thus we have 
$$
TV_0^{(N+1)\tau} (\dot\chi^h) = O(1)\sum_{m=1}^{N+1} |{\Cal E}(\Delta_{0,m})|,
\leqno{(4.7)}
$$
with
$$
\sum_{m=1}^{m(0)} |{\Cal E}(\Delta_{0,m})| = O(1) \, TV^0_{-(T/\lambda)-2h}
(v_L^0-c_1 w_L^0)
\leqno{(4.8a)}
$$
and
$$
\sum_{m=m(0)} ^{N+1} |{\Cal E}(\Delta_{0,m})| =\cases
O(1)\, TV_0^{(T/\lambda)+2h} (u_R^0)
\hfil & \text{ if } \quad V^*\ne c_3\cr
O(1) TV_0 ^{(T/\lambda)+2h}(u_R^0)+O(1)N
 \hfil & \text{ if } \quad V^*=c_3.\cr
\endcases 
\leqno{(4.8b)}
$$
Combining (4.7) and (4.8) gives (4.3). The proof of the lemma is complete.
~~ $\Box$

\vskip.2cm

\

%




\vskip.2cm

\noindent{\bf Proof of Theorem 4.1.} In view of Lemma~4.1 and Ascoli's theorem,
we have the convergence result (4.4). In view of Lemma~4.2, the total
variation of $\dot\chi^h$ on a compact set $[0,T]$ is uniformly bounded.
By extracting again a subsequence, Helly's Theorem gives
$$
\dot \chi^h(t)\to k(t)\quad \hbox{ for all times } t \ge 0
\leqno{(4.9a)}
$$
and
$$
TV_0^t (\dot\chi^h) \to \ell (t) \quad\hbox{ for all times } t\ge 0,
\leqno{(4.9b)}
$$
where $k\in BV_{loc}([0,\infty[,\RR)$ and $\ell \in L_{loc}^\infty([0,\infty
[,\RR_+)$. We define $E\subset \RR_+$ as the set of all points of
discontinuity of the function $\ell$. This set is at most countable since
the function $\ell$ is non-decreasing (and so has bounded variation).

We are going to prove that (4.4) holds with the above choice of set $E$.
Let $t$ be in $\RR_+\setminus E$, and let $\epsilon  > 0$ be so small that
$$
TV_{t-\epsilon}^{t+\epsilon} (\dot\chi^h) < \epsilon.
\leqno{(4.10)}
$$
This is possible because $t\not\in E$. Then, in view of (4.9a) and (4.10),
we have
$$
k(t)-\epsilon < \dot\chi^h(s) < k(t)+\epsilon \quad \hbox{ for }
s\in ]t-\epsilon,t+\epsilon[.
\leqno{(4.11)}
$$
On the other hand, we know that the curve $\chi^h$ has the slope
$\dot\chi^h(t)$ on the interval $[n\tau, (n+1)\tau[\ni t$ and jumps
by $\pm 2h$ at times $(n+1)\tau$. The slope $\dot\chi^h$ of $\chi^h$
is ``controlled'' by inequalities (4.11), while the jumps of $\chi^h$
are determined by the given sequence $\{a_n\}$.

Let $n'$ and $n''$ be two integers such that $(n'-1)\tau\le t' <n'\tau$
and $n''\tau \le t'' <(n''+1)\tau$, where $t-\epsilon< t' <t''<t+\epsilon$
are given. We set
$$
\Omega_+ =\{ m / m \hbox{ integer, } n'\le m \le n'' \hbox{ and }
a_m < (k(t)-\epsilon) {\tau\over h}\}
$$
and
$$
\Omega^* = \{ m | m \hbox{ integers, } n'\le m \le n'' \hbox{ and }
a_m > (k(t)+\epsilon) {\tau\over h}\}.
$$
In view of (4.11) and between times $n'\tau$ and $n''\tau$, the curve
$\chi^h$ has at least $\#\Omega_*$ jumps to the right and at most
$n''-n'-\#\Omega_*$ to the left, thus we have
$$
\chi^h (n''\tau) - \chi^h(n'\tau) \ge (2\# \Omega_* -n'+n'')h.
\leqno{(4.12a)}
$$
Similarly for $\Omega^*$ we get
$$
\chi^h(n''\tau) -\chi^h (n'\tau)\le (2\# \, \Omega_* -n'+n'')h.
\leqno{(4.12b)}
$$
But the equidistributedness of $\{a _n\}$ means that
$$
{\# \, \Omega_*\over n''-n'} \to {1\over 2} + (k(t)-\epsilon) {\tau
\over 2h} \quad\hbox{ and }\quad {\#\, \Omega^*\over n''-n'}\to {1\over 2}
(k(t)+\epsilon){\tau\over 2h}
\leqno{(4.13)}
$$
when $h\to 0$.

Combining (4.12) and (4.13) and letting $h\to 0$ yield the inequalities
$$
\chi (t'') -\chi (t') \ge (k(t)-\epsilon) \, (t''-t')
$$
and
$$
\chi (t'')-\chi (t')\le (k(t) + \epsilon )\, (t''-t'),
$$
which are valid for all $t-\epsilon < t' < t'' < t+\epsilon$
and thus in particular imply
$$
k(t)-\epsilon \le {d\chi (t')\over dt} \le k(t)+\epsilon
\quad\hbox{ for }\quad t-\epsilon < t' < t+\epsilon.
\leqno{(4.14)}
$$
Letting $\epsilon$ go to zero in (4.14) yields
$$
{d\chi (t)\over dt} = k(t).
$$
The proof is complete. ~~$\Box$

\

\noindent{\bf Remark 4.1.} Estimate (4.3) of Lemma~4.2 makes clear that only the
$c_1$ waves located at the left of the initial phase discontinuity and the 
$\pm c_3$ waves located at the right of the initial phase discontinuity 
contribute to the change in speed of the phase boundary.

\

We finally prove that the result in Theorem~4.1 is sufficient for the
passage to the limit in the kinetic relation.

\

\noindent{\bf Theorem 4.2.} {\it Suppose that} $V^*\ne 0$. {\it Then the limit function}
$u$ {\it given by Glimm's scheme satisfies}
$$
\partial _t U(u) +\partial _x F(u) = -\nu_t \phi(-{\nu_t\over\nu_x})
\delta _{x=\chi(t)}
\leqno{(4.15)}
$$
$H_1${\it --almost everywhere on the set } ${\Cal B}_{sub}(u)$.

\

Equality (4.15) is understood as equality between Borel measures on
$\RR_+\times \RR$. Here ${\Cal B}_{sub}(u)$ (according to the definition
of Section~1) is the set of all points of approximate jump of $u$
associated with subsonic phase discontinuities.
In view of the formula of Section~1, it is clear that, when $V^*\ne 0$,
(4.15) is equivalent to the formulation (1.11) of the kinetic 
relation. The case
$V^* =0$ could in principle be treated by the same technique
but this would require further analysis.

\

\noindent{\bf Remark 4.2.} 1)~~The pointwise convergence property of Glimm's 
scheme was already used in LeFloch-Liu [30] to derive an existence 
result for nonlinear hyperbolic systems in nonconservative form.

\noindent 2)~~If $V^*=0$,  $\dot\chi$ may vanish and then relation (4.15)
is not sufficient to uniquely characterize the solution
(e.g. of the Riemann problem).

\vskip.2cm

\noindent{\bf Proof of Theorem 4.2.} For all times $t\ge 0$, we introduce an approximate normal
$\nu^h(t) \break = (\nu_t^h(t), \, \nu_x^h(t)\, )$ by
$$
\nu _t^h(t) ^2 +\nu_x^h(t)^2 = 1, \qquad
\dot\chi ^h (t) = -{\nu_t^h(t)\over \nu_x^h(t)} \quad\hbox{ and }\quad
\nu_x^h(t) > 0.
$$
Similarly, from $\dot \chi(t)$, we define $\nu(t) = (\nu_t(t),\nu_x(t))$.
According to the notation of Section~1, we have in fact
$\nu(t) = \nu(t,\chi(t)\, )$. First of all, we claim that
$$
\nu_t^h\phi (\dot\chi^h) \delta_{x=\chi^h} \quad \to \quad
\nu_t\phi ({d\chi\over dt}) \delta_{x=\chi}
\leqno{(4.16)}
$$
in the weak-star topology of bounded Borel measures on $\RR_+\times\RR$.

Since $\dot \chi ^h$ satisfies (4.3) and the right hand side of (4.3) is
small by the assumption (3.1b), the function $\dot\chi ^h$ has small total
variation. When $V^*\ne 0$, we can ensure that $\dot\chi^h$ is
bounded away from zero uniformly with respect to $h$. In view of 
properties (1.8), the function $\phi$ is (at least) continuous in the
range of values taken by $\dot\chi^h$. This fact combined with the result
of convergence (4.5) gives
$$
\nu_t^h(t) \, \phi(\dot\chi^h(t)\, )\to \nu_t(t)\, \phi({d\chi\over dt}
(t)\, ) \quad\hbox{ for all } \quad t\in \RR_+ \setminus E,
\leqno{(4.17)}
$$
where $E$ is an at most countable set.  From (4.17) and the uniform 
convergence of $\chi^h$ to $\chi$ Cf.~(4.4), we deduce that
$$
\int \nu_t^h (t)\, \phi(\dot\chi^h(t)\, )\, \theta (\chi^h(t)\, \, dt
\to \int \nu_t (t)\, \phi({d\chi\over dt}(t)\, )\, \theta(\chi (t)\, )\, dt
$$
for each continuous function $\theta\colon\, \RR\to\RR$ with compact support.
This proves (4.16).

By construction, the approximate solutions $u^h$ satisfy the kinetic relation
$$
\partial _t U(u^h) +\partial_x F(u^h) = -\nu_t^h\phi(-{\nu_t^h\over \nu_x^h}
)\, \delta_{x=\chi^h}
\leqno{(4.18)}
$$
$H_1$--almost everywhere on the set ${\Cal B}_{sub}(u^h)$. We claim that
using (4.16) we can pass to the limit in (4.18) and obtain
$$
\partial _t U(u) +\partial _x F(u) = - \nu_t \, \phi(-{\nu_t\over\nu_x})\,
\delta_{x=\chi}
\leqno{(4.19)}
$$
$H_1$--almost everywhere on the set ${\Cal B}_{sub}(u)$.

The left hand sides of (4.18) and (4.19) are treated easily since they
have a (divergence-like) conservation form. In particular, we have
$$
\partial _t U(u^h) +\partial_x F(u^h)\to \partial_t U(u) +\partial_x
F(u)
\leqno{(4.20)}
$$
in the weak star topology of bounded Borel measures on $\RR_+\times \RR$.

In case $V^*<c_3$, (4.18) is satisfied on the whole space
$\RR_+\times \RR$ and so the desired result (4.19) is an immediate
consequence of (4.18), (4.16) and (4.20).

When $V^*>c_3$, nothing has to be proved since no kinetic relation
is imposed then.

The final case $V^* = c_3$ is treated as follows. We note that one can
find two Lipschitz continuous functions $\tilde\phi (V)$ and 
$\tilde \phi_+(V)$ defined for $V$
in a neighborhood of $c_3$ such that the kinetic relation (e.g. for
$u^h$) is equivalent to the two inequalities
$$
\partial_t U(u^h) + \partial _x F(u^h) \le -\nu_t^h \tilde\phi_+(-{\nu_t^h
\over \nu_x^h})\, \delta_{x=\chi^h}
\leqno{(4.21a)}
$$
and
$$
\partial _t U(u^h) + \partial_x F(u^h)\ge -\nu_t^h \tilde\phi_-(-{\nu_t^h
\over \nu_x^h})\, \delta_{x=\chi^h},
\leqno{(4.21b)}
$$
where $\tilde \phi_\pm$ are chosen in such a way that
$$
\tilde \phi_+(V) =\tilde\phi_-(V)=\phi(V) \quad\hbox{ for }\quad V<c_3
$$
and
$$
\tilde \phi_+ \hbox{ and } \tilde\phi_- \hbox{ are Lipschitz continuous with:}
\quad
\tilde \phi_-(V) < \tilde\phi_+(V).
$$
Namely, this is possible since (4.21a), (4.21b) when $V\le c_3$ give back
the kinetic relation; while for $V>c_3$ (4.21a), (4.21b) are trivially
satisfied provided that the entropy dissipation in the supersonic
case remains in the interval $[\tilde\phi_-(V), \tilde\phi_+(V)]$.
In this latter case, the entropy dissipation across the phase
boundary, say $\tilde\phi(V)$, is the following (cf.~the notation of
Section~2):
$$
\tilde\phi (V) ={1\over 2} (k_1-k_3) (w_Mw_m-w_Rw_-) = {1\over 2}
(k_1-k_3)(w_Mw_m -{c_3+V\over c_1+V} w_Rh_{LR}),
$$
where $V=V(h_{LR}) > c_3$ is a root of the equation (2.12b).
By (1.8b) and (1.9b), we have
$$
\lim_{u\to c_3\atop V>c_3} \tilde\phi(V) = {1\over 2} (k_1-k_3)
w_M w_m =\bar\psi (c_3) =\lim_{V\to c_3\atop V<c_3} \phi(V).
$$
This proves the continuity of the entropy dissipation at $V=c_3$.
Moreover $\tilde\phi$ is clearly Lipschitz continuous in view of Lemma~2.2.

Hence, for $\tilde\phi_\pm(V)$ suitably chosen and $V-c_3$
sufficiently small, the entropy dissipation $\tilde\phi(V)$ remains
in the interval $[\tilde\phi_-(V),\tilde\phi_+(V)]$.

It is clear that (4.16) still holds if $\phi$ is replaced by $\tilde \phi_-$
or $\tilde\phi_+$, i.e. we have in the weak-star topology:
$$
\nu_t^h\tilde\phi_\pm (\dot\chi ^h) \delta_{x=\chi^h}
{\buildrel {\text{weak}} ~ *\over \longrightarrow}
\nu_t\tilde\phi_\pm ({d\chi\over dt})\delta_{x=\chi}.
\leqno{(4.22)}
$$
Then (4.20) and (4.22) used in (4.21) yield:
$$
\partial _t U(u) +\partial _x F(u) \le -\nu_t \tilde\phi_+ (-{\nu_t
\over \nu_x}) \delta_{x=\chi^h}
$$
and
$$
\partial _t U(u) +\partial _x F(u) \ge -\nu_t \tilde\phi_- (-{\nu_t
\over \nu_x}) \delta_{x=\chi^h}
$$
which give (4.15). The proof is complete. ~~$\Box$

\vskip.2cm

We summarize in the following theorem  the results obtained along Section~3 
and in the present section.

\vskip.2cm

\noindent{\bf Theorem 4.3.} {\it Consider the mixed system} (1.1) {\it with an
initial condition which is a small perturbation in the BV norm of a single 
propagating phase boundary with speed} $V^*$. {\it Suppose that} 
$V^*\ne 0$ {\it and condition} (3.2b) {\it is satisfied if} $V^*=c_3$. 
{\it Then Glimm's scheme for this problem converges to an admissible 
weak solution which has the structure described in} (3.3), (3.4).

\vskip.2cm

\noindent{\bf Remark 4.3.} 1) ~
Note that the proof of Theorem~4.2 uses the property that the entropy
dissipation across a contact discontinuity is identically zero.

\noindent 2)~We believe that Theorem~4.3 could be extended to a
finite number of propagating phase boundaries. Also the restriction
$V^*\ne 0$ is only a technical assumption and could be removed by 
using other techniques from [16].

\noindent 3)~However, there is a main obstacle to a general result of
existence of BV solutions for (1.1). Indeed, for arbitrary large data,  
the phenomenon of initiation of new phases arises, and it is an open problem 
to derive a uniform bound on the total variation in that case. 

\vskip.2cm
\vfill

\

\centerline{\bf REFERENCES}
\frenchspacing
\vskip.4cm
\baselineskip.4cm
\item{[1]}
ABEYARATNE, R., KNOWLES, J. K., On the dissipative response due to
discontinuous strains in bars of unstable elastic materials, 
Int. J. Solids Structures 24(10), 1988, pp.~1021--1044.
\vskip.1cm
\item{[2]}
ABEYARATNE, R., KNOWLES, J. K., Kinetic relations and the propagation
of phase boundaries in solids, Arch. Rat. Mech. Anal. 114(2), 1991, 
pp.~119--154.
\vskip.1cm
\item{[3]}
ABEYARATNE, R., KNOWLES, J. K., Implications of viscosity and strain
gradient effects for the kinetics of propagating phase boundaries in
solids, SIAM J. Appl. Math. 51(5), 1991, pp.~1205--1221.
\vskip.1cm
\item{[4]}
ABEYARATNE, R., KNOWLES, J. K., On the propagation of maximally dissipative
phase boundaries in solids, Technical Report No.~1, California Inst.
of Technology, Pasa\-dena, CA, August 1990; to appear in Quart. of Appl. Math.
\vskip.1cm
\item{[5]}
CHERN, I. L., Stability theorem and truncation error analysis for the Glimm
scheme and for a front tracking method for flows with strong discontinuities,
Comm. Pure Appl. Math. 42, 1989, pp.~815--844.
\vskip.1cm
\item{[6]}
DAFERMOS, C. M., The entropy rate admissibility criterion for solutions
of hyperbolic conservation laws, Jour. Diff. Eq. 14, 1973, pp.~202--212.
\vskip.1cm
\item{[7]}
DAFERMOS, C. M., Hyperbolic systems of conservation laws, in ``Systems
of Nonlinear Partial Differential Equations'', ed. J. M. Ball, NATO
ASi Series C, No.~111, Dordrecht D. Reidel, 1983, pp.~25--70.
\vskip.1cm
\item{[8]}
DAL MASO, G., LeFLOCH, P., MURAT, F., Definition and weak stability
of non-conser\-vative products, preprint, Ecole Polytechnique (Palaiseau),
to be submitted.
\vskip.1cm
\item{[9]}
DI PERNA, R. J., Decay and asymptotic behavior of solutions to nonlinear
hyperbolic systems of conservation laws, Ind. Univ. Math. J. 24,
1975, pp.~1041--1071.
\vskip.1cm
\item{[10]}
DI PERNA, R., Singularities of solutions of nonlinear hyperbolic systems of
conservation laws, Arch. Rat. Mech. Anal. 60, 1975, pp.~75--100.
\vskip.1cm
\item{[11]}
DI PERNA, R. J., MAJDA, A., The validity of nonlinear geometrical optics
for weak solutions of conservation laws, Comm. Math. Phys. 98, 1985,
pp.~313--347.
\vskip.1cm
\item{[12]}
FAN, H., A vanishing viscosity approach on the dynamic of phase transitions
 in van der Waals fluids, IMA Preprint, 1991.
\vskip.1cm
\item{[13]}
FAN, H., SLEMROD, M., The Riemann problem for systems of conservation laws 
of mixed type, IMA Preprint, 1991.
\vskip.1cm
\item{[14]}
FEDERER, H., Geometric Measure Theory, Springer-Verlag, New York, 1969.
\vskip.1cm
\item{[15]}
GLIMM, J., Solutions in the large for nonlinear hyperbolic systems of
equations, Comm. Pure Appl. Math. 18, 1965, pp.~697--715. 
\item{[16]}
GLIMM, J., LAX, P. D., Decay of solutions of systems of nonlinear hyperbolic
conservation laws, Memoirs of the A.M.S., No.~101, 1970.
\vskip.1cm
\item{[17]}
GURTIN, M.E., On a theory of phase transitions with interfacial energy, 
Arch. Rat. Mech. Anal. 87 (1984), pp.~187--212.
\vskip.1cm
\item{[18]}
HAGAN, R., SLEMROD, M., The viscosity-capillarity criterion for shocks
and phase transitions, Arch. Rat. Mech. Anal. 83, 1983, pp.~333--361.
\vskip.1cm
\item{[19]}
HATTORI, H., The Riemann problem for a van der Waals fluid with entropy
rate admissibility criterion: Isothermal case, Arch. Rat. Mech. Anal. 92,
1986, pp.~246--263.
\vskip.1cm
\item{[20]}
HATTORI, H., The Riemann problem for a van der Waals fluid with entropy
rate admissibility criterion: Non-isothermal case, Jour. Diff. Eq. 65,
1986, pp.~158--174.
\vskip.1cm
\item{[21]}
HOU, L. Y., LeFLOCH, P., Why non-conservative schemes converge to wrong
solutions: error analysis, Math. of Comp., 1993.
\vskip.1cm
\item{[22]}
HSIAO, L., Uniqueness of admissible solutions of the Riemann problem for a 
system of conservation laws of mixed type, Jour. Diff. Eq. 86, 1990, 
pp.~197--233.
\vskip.1cm
\item{[23]}
JAMES, R. D., The propagation of phase boundaries in elastic bars,
Arch. Rat. Mech. Anal. 73, 1980, pp.~125--158.
\vskip.1cm
\item{[24]}
KEYFITZ, B. L., The Riemann problem for non-monotone stress-strain functions:
a ``hysteresis'' approach, Lec. in Appl. Math., Vol.~23, 1986, pp.~379--395.
\vskip.1cm
\item{[25]}
LAX, P. D., Hyperbolic systems of conservation laws, II, Comm. Pure Appl. Math. 10, 1957, pp.~ 537--566.
\vskip.1cm
\item{[26]}
LAX, P. D., Hyperbolic systems of conservation laws and the mathematical
theory of shock waves, Regional Conf. Series in Appl. Math., No.~11,
SIAM, Philadelphia, 1973.
\vskip.1cm
\item{[27]}
LeFLOCH, P.G., Entropy weak solutions to nonlinear hyperbolic systems in
non-conser\-vative form, Comm. Part. Diff. Eqs. 13(6), 1988, pp.~669--727.
\vskip.1cm
\item{[28]}
LeFLOCH, P.G., Shock waves for nonlinear hyperbolic systems in non-conservative
form, Preprint No.~599, Inst. of Math. and its Appli., Minneapolis, October
1989.
\vskip.1cm
\item{[29]}
LeFLOCH, P.G., Entropy weak solutions to nonlinear hyperbolic systems in conservative form,
 in Proc. Second Intern. Conf. on Nonlinear Hyperbolic problems,
Aachen, FRG, 1988, J. Ballmann and R. Jeltsch (Eds.),
 Note on Num. Fluid Mech. vol. 24, Vieweg, Braunschweig, 1989,
 pp.~362--373.
\vskip.1cm
\item{[30]}
LeFLOCH, P.G., LIU, T. P., Existence theory for nonlinear hyperbolic systems
in non-conservative form, Forum Mathematicum, 1992.
\vskip.1cm
\item{[31]}
LeFLOCH, P.G., XIN, Z. P., Uniqueness via the adjoint problem for some systems of
 conservation laws, Preprint Courant Institute, New York University; Comm. Pure Appl. Math., 1992.
\vskip.1cm
\item{[32]}
LIU, T. P., The Riemann problem for general systems of conservation laws, J. Diff. Equa. 18, 1975, pp.~218--234.
\vskip.1cm
\item{[33]}
LIU, T. P., The deterministic version of the Glimm scheme, Comm. Math. Phys. 57, 1977, pp.~135--148.
\vskip.1cm
\item{[34]}
LIU, T. P., Admissible solutions of hyperbolic systems of conservation laws, 
Memoirs of the AMS No.~240, 1981, pp~1--78.
\vskip.1cm
\item{[35]}
PEGO, R. L., Phase transitions in one-dimensional nonlinear viscoelasticity:
admissibility and stability, Arch. Rat. Mech. Anal. 97(4), 1987, pp.~353--394.
\vskip.1cm
\item{[36]}
PEGO, R. L., SERRE, D., Instability in Glimm's scheme for two systems
of mixed type, SIAM J. Num. Anal. 25, 1988, pp.~965--988.
\vskip.1cm
\item{[37]}
SCHOCHET, S., Sufficient conditions for local existence via Glimm's
scheme for large BV data, J. Diff. Equa. 89(2), 1991, pp.~317--354.
\vskip.1cm
\item{[38]}
SHEARER, M., The Riemann problem for a class of conservation laws of mixed
type, Jour. Diff. Eq. 46, 1982, pp.~426--443.
\vskip.1cm
\item{[39]}
SHEARER, M., Nonuniqueness of admissible solutions of the Riemann
initial value problem for a system of conservation laws of mixed type,
Arch. Rat. Mech. Anal. 93, 1986, pp.~45--59.
\vskip.1cm
\item{[40]}
SHEARER, M., Dynamic phase transitions in a van der Waals gas, Quart. Appl.
Math. 46(4), 1988, pp.~631--636.
\vskip.1cm
\item{[41]}
SLEMROD, M., Admissibility criteria for propagating phase boundaries in
a van der Waals fluid, Arch. Rat. Mech. Anal. 81, 1983, pp.~301--315.
\vskip.1cm
\item{[42]}
SLEMROD, M., A limiting viscosity approach to the Riemann problem for
materials exhibiting change of phase, Arch. Rat. Mech. Anal. 105, 1989,
pp.~327--365.
\vskip.1cm
\item{[43]}
TEMPLE, B., Global solution of the Cauchy problem for a class of $2\times 2$
non-strictly hyperbolic conservation laws, Adv. in Appl. Math. 3, 1982,
pp.~335--375.
\vskip.1cm
\item{[44]}
TRUSKINOVSKY, L., Kinks versus shocks, in ``Shock induces transitions
and phase structures in general media'', ed. R. Fosdick, E. Dunn and
H. Slemrod, Springer-Verlag, 1991.
\vskip.1cm
\item{[45]}
VOLPERT, A. I., The space BV and quasilinear equations, USSR Sb.~2,
1967, pp.~225--267.
\vskip.1cm

%


\vfill
\eject
\end